\documentclass[12pt]{article}
\usepackage[dvips]{epsfig}
\usepackage{amsthm,amsmath,amssymb}
\usepackage{latexsym}
\usepackage{placeins}
\usepackage{color}
\usepackage[colorlinks]{hyperref}
\definecolor{blue}{rgb}{0,0,0.9}
\definecolor{red}{rgb}{0.9,0,0}
\definecolor{green}{rgb}{0,0.50,0.10}
\definecolor{violet}{rgb}{0.5804,0.0000,0.8275}

\newcommand{\red}[1]{\begin{color}{red}#1\end{color}}

\def\@themcountersep{}

\newtheorem{THEO}{Theorem}[section]
\newtheorem{ALGo}[THEO]{Algorithm}
\newtheorem{CONJ}[THEO]{Conjecture}
\newtheorem{COND}[THEO]{Condition}
\newtheorem{ASSUMP}[THEO]{Assumption}
\newtheorem{CORO}[THEO]{Corollary}
\newtheorem{DEFI}[THEO]{Definition}
\newtheorem{EXAMP}[THEO]{Example}
\newtheorem{FACT}[THEO]{Fact}
\newtheorem{HYPO}[THEO]{Hypothesis}
\newtheorem{LEMM}[THEO]{Lemma}
\newtheorem{PROB}[THEO]{Problem}
\newtheorem{PROP}[THEO]{Proposition}
\newtheorem{REMA}[THEO]{Remark}
\newcommand{\theo}{\begin{THEO}}
\newcommand{\algo}{\begin{ALGo} \rm}
\newcommand{\cond}{\begin{COND} \rm}
\newcommand{\assump}{\begin{ASSUMP} \rm}
\newcommand{\conj}{\begin{CONJ}}
\newcommand{\coro}{\begin{CORO}}
\newcommand{\defi}{\begin{DEFI} \rm}
\newcommand{\examp}{\begin{EXAMP} \rm}
\newcommand{\fact}{\begin{FACT}}
\newcommand{\hypo}{\begin{HYPO} \rm}
\newcommand{\lemm}{\begin{LEMM}}
\newcommand{\prob}{\begin{PROB} \rm}
\newcommand{\prop}{\begin{PROP}}
\newcommand{\rema}{\begin{REMA} \rm}
\newcommand{\etheo}{\end{THEO}}
\newcommand{\ealgo}{\end{ALGo}}
\newcommand{\econd}{\end{COND}}
\newcommand{\eassump}{\end{ASSUMP}}
\newcommand{\econj}{\end{CONJ}}
\newcommand{\ecoro}{\end{CORO}}
\newcommand{\edefi}{\end{DEFI}}
\newcommand{\eexamp}{\end{EXAMP}}
\newcommand{\efact}{\end{FACT}}
\newcommand{\ehypo}{\end{HYPO}}
\newcommand{\elemm}{\end{LEMM}}
\newcommand{\eprob}{\end{PROB}}
\newcommand{\eprop}{\end{PROP}}
\newcommand{\erema}{\end{REMA}}

\def\0{\mbox{\bf 0}}
\def\1{\mbox{\bf 1}}
\def\2{\mbox{\bf 2}}
\def\3{\mbox{\bf 3}}
\def\4{\mbox{\bf 4}}
\def\5{\mbox{\bf 5}}
\def\6{\mbox{\bf 6}}
\def\7{\mbox{\bf 7}}
\def\8{\mbox{\bf 8}}
\def\9{\mbox{\bf 9}}

\def\b{\mbox{\boldmath $b$}}
\def\cc{\mbox{\boldmath $c$}}
\def\d{\mbox{\boldmath $d$}}

\def\u{\mbox{\boldmath $u$}}

\def\w{\mbox{\boldmath $w$}}
\def\x{\mbox{\boldmath $x$}}

\def\A{\mbox{\boldmath $A$}}
\def\B{\mbox{\boldmath $B$}}
\def\C{\mbox{\boldmath $C$}}

\def\H{\mbox{\boldmath $H$}}

\def\O{\mbox{\boldmath $O$}}
\def\P{\mbox{\boldmath $P$}}
\def\Q{\mbox{\boldmath $Q$}}

\def\X{\mbox{\boldmath $X$}}
\def\Y{\mbox{\boldmath $Y$}}
\def\Z{\mbox{\boldmath $Z$}}
\def\AC{\mbox{$\cal A$}}
\def\BC{\mbox{$\cal B$}}

\def\TC{\mbox{$\cal T$}}

\def\inprod#1#2{\langle#1, \, #2\rangle}

\def\Real{\mbox{$\mathbb{R}$}}
\def\coneK{\mbox{$\mathbb{K}$}}
\def\coneJ{\mbox{$\mathbb{J}$}}

\def\spaceV{\mbox{$\mathbb{V}$}}
\def\spaceL{\mbox{$\mathbb{L}$}}
\def\SymMat{\mbox{$\mathbb{S}$}}

\def\SymN{\mbox{$\mathbb{N}$}}
\def\Integer{\mbox{$\mathbb{Z}$}}

\def\CPP{\mbox{$\mathbb{CPP}$}}
\def\CP{\mbox{$\mathbb{COP}$}}

\def\balpha{\mbox{\boldmath $\alpha$}}
\def\bbeta{\mbox{\boldmath $\beta$}}
\def\salpha{\mbox{\scriptsize $\balpha$}}
\def\sbeta{\mbox{\scriptsize $\bbeta$}}
\def\bgamma{\mbox{\boldmath $\gamma$}}
\def\sgamma{\mbox{\scriptsize $\bgamma$}}

\def\bdelta{\mbox{\boldmath $\delta$}}
\def\sdelta{\mbox{\scriptsize $\bdelta$}}

\def\bGamma{\mbox{$\bf{\Gamma}$}}

\def\sAC{\mbox{\scriptsize $\AC$}}
\def\sBC{\mbox{\scriptsize $\BC$}}

\def\inprod#1#2{\langle #1,\,#2\rangle}

\begin{document}

\title{ \large 
A Geometrical Analysis of a Class of Nonconvex Conic Programs for Convex Conic Reformulations of Quadratic and Polynomial Optimization Problems 
}

\bigskip
\author{
\normalsize 
Sunyoung Kim\thanks{Department of Mathematics, Ewha W. University, 52 Ewhayeodae-gil, Sudaemoon-gu, Seoul 03760, Korea 
			({\tt skim@ewha.ac.kr}). The research was supported
               by   2017-R1A2B2005119.}, \and \normalsize
Masakazu Kojima\thanks{Department of Industrial and Systems Engineering,
	Chuo University, Tokyo 192-0393, Japan 
({\tt kojima@is.titech.ac.jp}).
	This research was supported by Grant-in-Aid for Scientific Research (A) 26242027.},
%	and the Japan Science and 
%	Technology Agency (JST), the Core Research of Evolutionary Science and 
%	Technology (CREST) research project.             
% 	},
  \and \normalsize
Kim-Chuan Toh\thanks{Department of Mathematics, and Institute of Operations Research and Analytics, National University of Singapore,
10 Lower Kent Ridge Road, Singapore 119076
({\tt mattohkc@nus.edu.sg}). 
This research is supported in part by the Ministry of Education, Singapore, Academic Research Fund (Grant number: R-146-000-257-112).
 } 
}
\date{\normalsize \today}

\maketitle 
\vspace*{-0.8cm}

\begin{abstract}
\noindent
We present a geometrical analysis on the completely positive programming reformulation of quadratic 
optimization problems and its extension to polynomial optimization problems with a class of geometrically defined nonconvex 
conic programs and their covexification. The class of nonconvex conic programs is described with a linear objective function
in a linear space $\spaceV$, and the constraint set  is represented  
geometrically as the intersection of a nonconvex cone $\coneK \subset \spaceV$, a face $\coneJ$ of the convex hull of $\coneK$  
and a parallel translation $\spaceL$ of a supporting hyperplane of the nonconvex cone $\coneK$.  
We show that under a moderate assumption, the original nonconvex conic program can equivalently be reformulated as a convex
conic program  by replacing the constraint set with the intersection of $\coneJ$ and 
the hyperplane $\spaceL$. The replacement procedure is  applied to derive the completely positive programming reformulation of quadratic 
optimization problems and its extension to polynomial optimization problems. 
\end{abstract}


\vspace{0.5cm}


%!TEX root = ./main.tex
\section{Introduction}

Polynomial optimization problems (POPs) is a major
class of optimization problems in theory and practice. 
Quadratic optimizations problems (QOPs) are, in particular, a widely studied subclass of POPs 
as they include many important NP-hard combinatorial  problems such as 
binary QOPs, maximum 
stable set problems, graph partitioning problems and quadratic assignment problems. 
To numerically solve QOPs, a common approach is through solving their convex conic relaxations such as semidefinite programming relaxations \cite{POLJAK1995,LASSERRE2001} 
and doubly nonnegative (DNN) relaxations \cite{GE2010,KIM2013,TNW2012,YOSHISE2010}. 
As those relaxations provide lower bounds of different qualities, 
the tightness of the lower bounds 
has been a very critical issue in assessing the strength of the relaxations.
The completely positive programming (CPP) reformulation of QOPs, 
which provides their exact optimal values, has been extensively studied in theory.
More specifically, 
QOPs over the standard simplex \cite{BOMZE2000,BOMZE2002}, 
maximum stable set problems \cite{DEKLERK2002}, 
graph partitioning problems \cite{POVH2007}, 
and quadratic assignment problems \cite{PRendl09} are equivalently  reformulated as CPPs. 
Burer's CPP reformulations \cite{BURER2009} of a class of linearly constrained QOPs 
in nonnegative and binary variables 
provided a more general framework to study  the specific problems mentioned above. See also the papers 
\cite{ARIMA2012,ARIMA2013,BOMZE2017,DICKINSON2012,PENA2015} for further developments.

Despite a great deal of studies on the CPP relaxation,
its geometrical aspects have not been well understood. 
The main purpose of this paper is to present and analyze  {\em essential features of the CPP reformulation of QOPs and its extension to POPs 
by investigating their geometry.} 
With the geometrical analysis, many existing equivalent  reformulations of QOPs and POPs can be considered in a 
unified manner 
 and  deriving effective numerical methods for computing tight bounds can be facilitated.
In particular, the class of QOPs that can be equivalently reformulated as CPPs in our framework includes
Burer's  
class of linearly constrained QOPs in nonnegative and binary variables \cite{BURER2009}  as a special case;
see Sections 2.2 and 6.1.
% \blue{Our work also subsume the analysis done in \cite{PENA2015} for POPs.}

%%%%%%%%%%%%%%%%%%%%%%
\subsection{A geometric framework for the CPP relaxation of QOPs and its extension to POPs}

A nonconvex conic optimization problem (COP), denoted as COP($\coneK_0,\Q^0$), of the form presented below 
is the most distinctive feature of our framework for the CPP relaxation of QOPs and its extension to POPs. 
Let $\spaceV$ be a finite dimensional vector space with the inner product $\inprod{\A}{\B}$ for every pair of $\A$ and $\B$ 
in $\spaceV$. 
For a cone $\coneK \subset \spaceV$,
let co$\coneK$ denote the convex hull of $\coneK$ 
and $\coneK^*$ the dual of $\coneK$, {\it i.e.}, $\coneK^* = \left\{\Y \in\spaceV: \inprod{\X}{\Y} \geq 0 \ \mbox{for every } \x \in \coneK \right\}$. 
Let $\H^0 \in \spaceV$, which will be described more precisely in Section 2.2 for QOPs and in Section 5 for general POPs. 
For every cone $\coneK_0 \subset \spaceV$ (not necessarily convex nor closed) and $\Q^0 \in \spaceV$, 
we consider the COP given by
\begin{eqnarray*}
\mbox{COP($\coneK_0,\Q^0$): } \zeta = \inf\left\{\inprod{\Q^0}{\X} : \X \in \coneK_0, \ \inprod{\H^0}{\X} = 1 \right\}. 
\end{eqnarray*}
Although this problem takes a very simple form, it plays a fundamental role throughout.
A key property is that COP($\coneK_0,\Q^0$) is equivalent to its {\it covexification}, 
COP(co$\coneK_0,\Q^0$) under the following conditions (Theorem~\ref{theorem:main0}). \vspace{2mm} \\ 
{\bf Condition I$_0$: }
COP($\coneK_0,\Q^0$) is feasible and 
 $\O \not= \H^0 \in \coneK_0^*$. \vspace{2mm}\\
{\bf Condition II$_0$: }
$\inf\left\{\inprod{\Q^0}{\X} : \X \in \coneK_0, \ \inprod{\H^0}{\X} = 0 \right\} \geq 0$. \vspace{2mm}\\
The only restrictive and essential condition among the conditions is $\O \not= \H^0 \in \coneK_0^*$, while the others are natural. 
It means that $\left\{ \X \in \spaceV : \inprod{\H^0}{\X} = 0 \right\}$ forms a 
supporting hyperplane of the cone $\coneK_0$ and that the feasible region of COP($\coneK_0,\Q_0$) is described as 
the intersection of the nonconvex cone $\coneK_0 \subset \spaceV$ and a parallel 
translation of the supporting hyperplane of $\coneK_0$. Condition II$_0$ is necessary to ensure that the optimal value of COP(co$\coneK_0,\Q^0$) is finite. See Figure 1 in Section 3.1 for illustrative examples 
of COP($\coneK_0,\Q^0$) which satisfies Condition I$_0$ and 
II$_0$.

We consider a specific COP($\coneK_0,\Q^0$) with 
$\coneK_0 = \coneK \cap \coneJ$ for some cone $\coneK \subset \spaceV$ and some face $\coneJ$ of co$\coneK$. 
Note that $\coneK\cap\coneJ$ is a nonconvex cone. 
Since $\coneJ$ is a face of co$\coneK$, we have that co$(\coneK\cap\coneJ) = \coneJ$ ((i) of Lemma~\ref{lemma:basic}).  
It follows that the equivalence of COP(co$(\coneK\cap\coneJ),\Q^0$) and COP($\coneJ,\Q^0$) holds trivially. This is another distinctive feature of our geometric framework. 

In this paper, we mainly deal with a class of general POPs of the form: 
\begin{eqnarray}
\zeta^* = \inf \left\{ f_0(\w) : \w \in \Real^n_+, \ f_i(\w) = 0 \ (i=1,\ldots,m) \right\}, \label{eq:POP10} 
\end{eqnarray}
where $\Real^n_+$ denotes the nonnegative orthant of the $n$-dimensional Euclidean space $\Real^n$ and 
$f_i(\w)$ a real valued polynomial function in $\w = (w_1,\ldots,w_n) \in \Real^n$ $(i=0,\ldots,m)$. 
When all $f_i(\w)$ $(i=0,\ldots,m)$ are quadratic functions, \eqref{eq:POP10} becomes %provides 
a class of QOPs 
considered in this paper.

The equivalence of  COP(co$(\coneK\cap\coneJ),\Q^0$) and COP($\coneJ,\Q^0$) shown above can be applied to POP~\eqref{eq:POP10}
by just reducing  POP~\eqref{eq:POP10} to the form of COP($\coneK\cap\coneJ,\Q^0$). 
This reduction is demonstrated in Section 2.2 for QOP cases, 
and in Section 5 for general POP cases. 
For the resulting COP(co$(\coneK\cap\coneJ),\Q^0$) to satisfy Conditions I$_0$ and II$_0$ with $\coneK_0 = \coneK\cap\coneJ$, 
certain assumptions must be imposed. 
For example, if the feasible region of POP~\eqref{eq:POP10} is nonempty and bounded, and 
$f_i(\w)$ $(i=1,\ldots,m)$ are nonnegative for every $\w \geq \0$, COP($\coneK\cap\coneJ,\Q^0$) can be constructed 
such that Conditions I$_0$ and II$_0$ are satisfied with $\coneK_0=\coneK\cap\coneJ$ for some cone $\coneK$ and some face 
$\coneJ$ of co$\coneK$. Consequently, COP($\coneJ,\Q^0$) 
is indeed a convex COP reformulation of POP~\eqref{eq:POP10} with the 
same objective value $\zeta = \zeta^*$ (Theorem~\ref{theorem:main50}). 
Note that co$\coneK$ 
corresponds the CPP cone when POP~\eqref{eq:POP10} is a QOP, while it corresponds to an extension of the CPP cone 
for a general
POP. Thus, $\coneJ$ is a face of the CPP cone in the QOP case or a face of the extended CPP cone in the 
general POP case.

In the convexification from POP~\eqref{eq:POP10} to COP($\coneJ,\Q^0$), the objective function $f_0(\w)$ is 
relaxed to the linear function $\inprod{\Q^0}{\X}$ in $\X \in \mbox{co}\coneK$. The problem
COP($\coneJ,\Q^0$), however, does not explicitly involve 
any linear equality in $\X \in \mbox{co}\coneK$  induced from each equality constraint $f_i(\w) = 0$ $(i=1,\ldots,m)$. 
In fact, the feasible region of COP($\coneK\cap\coneJ,\Q^0$) is  geometrically represented in terms of 
a nonconvex cone $\coneK$, a face $\coneJ$ of co$\coneK$ and a hyperplane $\left\{ \X \in \spaceV : \inprod{\H^0}{\X} = 1 \right\}$.
This formulation is essential 
to derive the convex COP reformulation COP($\coneJ,\Q^0$) of QOPs and POPs  in a simple geometric setting. 

%%%%%%%%%%%%%%%%%%%
\subsection{Relations to existing works}

The geometric framework mentioned in the previous section 
generalizes the authors' previous work 
\cite{ARIMA2012,ARIMA2013,ARIMA2018,ARIMA2018b,KIM2013}. 
A convex reformulation  of a nonconvex COP in a vector space $\spaceV$ was also discussed
and  the results obtained there were applied to QOPs in 
\cite{ARIMA2012,ARIMA2018,KIM2013}  and 
POPs in \cite{ARIMA2013,ARIMA2018b}. 
Unlike the current framework, 
a fundamental difference in the previous framework lies in utilizing  a nonconvex COP of the form 
\begin{eqnarray}
 \zeta = \inf\left\{
 \inprod{\Q^0}{\X} : 
\begin{array}{l}
\X \in \coneK, \ \inprod{\H^0}{\X} = 1, \\
 \inprod{\Q^p}{\X} = 0 \ (p=1,\ldots,m) 
 \end{array} 
 \right\}, \label{eq:COP11} 
\end{eqnarray}
where $\coneK \subset \spaceV$ denotes a cone, $\Q^p \in \spaceV$ $(p=0,\ldots,m)$ and $\H^0 \in \spaceV$. 
In \cite{ARIMA2018,ARIMA2018b,KIM2013}, they imposed the assumption that $\Q^p \in \coneK^*$ 
$(p=0,\ldots,m)$ in addition to $\O \not= \H^0 \in \coneK^*$ and a  condition similar
to %(or essentially the same condition as) 
Condition II$_0$. Under this assumption, 
\begin{eqnarray}
\coneJ = \left\{ \X \in \mbox{co}\coneK : 
 \inprod{\Q^p}{\X} = 0 \ (p=1,\ldots,m) \right\} \label{eq:coneJ10} 
\end{eqnarray}  
forms a face of co$\coneK$ (Lemma 2.1). However, the converse is not true. A face 
$\coneJ$ of co$\coneK$ can be represented as in~\eqref{eq:coneJ10}  by some 
$\Q^p \in \coneK^*$ $(p=1,\ldots,m)$ iff it is an exposed face of co$\coneK$; 
hence if $\coneJ$ is a non-exposed face of co$\coneK$, such a representation in terms of some 
$\Q^p \in \coneK^*$ $(p=1,\ldots,m)$ is impossible. Very recently, Zhang \cite{ZHANG2018} showed that the CPP cone with 
dimension not less than $5$ is not facially exposed, i.e., some of its faces are non-exposed (see also 
\cite{BERMAN2015,DICKINSON2011} for geometric properties of the CPP cone). 
Thus, our framework using 
COP(co$(\coneK\cap\coneJ),\Q^0$) is more general than the work using 
\eqref{eq:COP11}  in  % \cite{ARIMA2018,ARIMA2018b,KIM2013}.}
\cite{ARIMA2012,ARIMA2013,ARIMA2018,ARIMA2018b,KIM2013}. 

The class of QOPs that can be reformulated as equivalent CPPs of 
the form COP($\coneJ,\Q^0$) in our framework covers most 
of the known classes of QOPs that can be reformulated as CPPs mentioned above, 
including Burer's class 
\cite{BURER2009} of linearly constrained QOPs in nonnegative and binary variables. 
With respect to  extensions to POPs presented in \cite{ARIMA2013,ARIMA2018b,PENA2015},
our geometric framework using  COP($(\coneK\cap\coneJ,\Q^0$) can be regarded as a generalization of the framework proposed in 
\cite{ARIMA2013,ARIMA2018b} where a class of POPs of the form~\eqref{eq:POP10}  is reduced to COP~\eqref{eq:COP11}. 
In~\cite{PENA2015}, Pe$\tilde{{\rm n}}$a, Vera and Zuluaga introduced the cone of completely positive tensor 
as an extension of the completely positive matrix for 
deriving equivalent convex relaxation of POPs.
 The class of POPs that can be convexified 
using their completely positive tensor cone is similar to our class that can be reformulated as equivalent CPPs of 
the form COP($\coneJ,\Q^0$). In fact, one of the two conditions imposed on their class, (i) of Theorem~4 in~\cite{PENA2015},  
corresponds to our condition~\eqref{eq:CondPOP1-1}, which was originated from a hierarchy of copositivity condition 
proposed in
\cite{ARIMA2012}. The other condition using ``the horizon cone'' in  (ii) of Theorem~4 of~\cite{PENA2015}, is different 
from our condition~\eqref{eq:CondPOP1-2}, but they are similar in nature (see Section~6 of \cite{ARIMA2012}). 
We should mention, however, that our framework is quite different form theirs. 

The above discussions show the versatility of our geometric framework in that it is applicable
to almost all known equivalent reformulations of QOPs as well as the more general case of POPs.

%%%%%%%%%%%%%%%%%%%%%%%
\subsection{Outline of the paper}

After introducing some notation and symbols in Section 2.1, we present how 
a general QOP can be reduced to COP($\coneK_0,\Q^0$) in Section 2.2, and present some 
fundamental properties of cones and their faces in Section 2.3. We establish the equivalence of COP($\coneK_0,\Q_0$) and its 
convexification COP(co$\coneK_0,\Q^0$) under Conditions I$_0$ and II$_0$ in Section 3.1, 
and derive the equivalence of COP($\coneK\cap\coneJ,\Q_0$) and its 
convexification COP($\coneJ,\Q^0$) by taking $\coneK_0 = \coneK \cap \coneJ$ for some cone $\coneK \subset \spaceV$ 
and some face $\coneJ$ of co$\coneK$ 
in Section 3.2.  In Section 4.1, we introduce 
a hierarchy of copositivity condition to represent a face $\coneJ$ of the convex hull 
co$\coneK$ of a cone $\coneK \subset \spaceV$ as in~\eqref{eq:coneJ10}. 
This  connects two  forms of a nonconvex COP,  COP($\coneK\cap\coneJ,\Q^0$) and COP~\eqref{eq:COP11}. 
In Section 4.2, some sufficient conditions for $\Q^p$ $(p=1,\ldots,m)$ to represent a face $\coneJ$ of co$\coneK$ as in~\eqref{eq:coneJ10}
are provided. Section 5 discusses convex COP reformulations of 
POPs as applications of the results obtained in Sections 3 and 4. We discuss homogenizing polynomials 
and an extension of the completely positive cone in Sections 5.1 and 5.2, respectively. We then construct a convex COP reformulation 
of POP~\eqref{eq:POP10} in Sections 5.3 and 5.4.  
In Section 6, we illustrate how we can apply the main theorems established in Section 5.4 to QOPs and POPs through 
examples. Finally, we conclude the paper in Section 7.

%!TEX root = ./main.tex
\section{Preliminaries}

\subsection{Notation and symbols}

Let $\Real^n$ denote the $n$-dimensional Euclidean space consisting of 
column vectors $\w = (w_1,\ldots,w_n)$, $\Real^n_+$ the nonnegative orthant of $\Real^n$, 
$\SymMat^n$ the linear space of $n \times n$ symmetric 
matrices with 
the inner product $\inprod{\A}{\B} = \sum_{i=1}^n\sum_{j=1}^n A_{ij} B_{ij}$, and $\SymMat^n_+$ the cone of positive semidefinite matrices in $\SymMat^n$. 
$\Integer^n$
denotes the set of integer column vectors in $\Real^n$, and 
$\Integer^n_+ = \Real^n_+ \cap \Integer^n$. 
$\cc^T$ denotes the transposition of 
a column vector $\cc \in \Real^n$. When $\Real^{1+n}$ is used, the first coordinate of $\Real^{1+n}$ is indexed by $0$ and 
$\x \in \Real^{1+n}$ is written as $\x = (x_0,x_1,\ldots,x_n) =  (x_0,\w) \in \Real^{1+n}$ with $\w \in \Real^{n}$. Also each matrix $\X \in \SymMat^{1+n} 
\subset \Real^{1+n} \times \Real^{1+n}$ has elements $X_{ij}$ $(i=0,\ldots,n, j=0,\ldots,n)$.

Let $\spaceV$ be a finite dimensional linear space with the inner product $\inprod{\A}{\B}$ for every pair of $\A$ and $\B$ in $\spaceV$ 
and $\left\|\A\right\| = \inprod{\A}{\A}^{1/2}$ for every $\A$ in $\spaceV$.
We say that $\coneK \subset \spaceV$ is a {\it cone}, 
which is not necessarily convex nor closed, 
if $\lambda\A \in \coneK$ for every $\A \in \coneK$ and $\lambda \geq 0$. 
Let  co$\coneK$ denote the convex hull of a cone $\coneK$, and 
cl$\coneK$ the closure of $\coneK$. $\SymMat^{1+n}$ may be regarded as a special case of a linear space $\spaceV$ in the 
subsequent discussions. 
Since $\coneK$ is a cone, we see that 
co$\coneK = \left\{\sum_{p=1}^m \X^p :  \X^p \in \coneK \ (p=1,\ldots,m) \ \mbox{for some } m \in\Integer_+ \right\}$. 
The 
{\it dual} of a cone $\coneK$ is defined as 
$\coneK^* = \left\{\Y \in \spaceV : \inprod{\Y}{\X} \geq 0 \ \mbox{for every } \ \X \in \coneK\right\}$. From the definition, 
we know that $\coneK^* = (\mbox{co}\coneK)^*$. It is well-known and also easily proved by the separation theorem of convex sets 
that $\coneK^{**} = \mbox{cl co}\coneK$, the closure of ${\rm co}\coneK.$

We note that a cone $\coneK$ is convex iff 
$\X = \sum_{i=1}^m \X^i \in \coneK$ whenever $\X^i \in \coneK$ $(i=1,\ldots,m)$. 
Let $\coneK$ be a convex cone in a linear space $\spaceV$. 
A convex cone $\coneJ \subset \coneK$ is said to be a {\it face} of $\coneK$ 
if $\X^1 \in \coneJ$ and $\X^2 \in \coneJ$ whenever 
$\X = \X^1/2 + \X^2/2 \in \coneJ$, $\X^1 \in \coneK$ and $\X^2 \in \coneK$ 
(the standard definition of a face of a convex set), or, 
if $\X^i \in \coneJ$ $(i=1,\ldots,m)$ whenever $\X = \sum_{i=1}^m \X^i \in \coneJ$ and 
$\X^i \in \coneK$ $(i=1,\ldots,m)$ (the equivalent characterization of a face of a convex cone). The equivalence can be easily shown by induction. 
A face $\coneJ$ of $\coneK$ is {\it proper} if 
$\coneJ \not= \coneK$, and a proper face $\coneJ$ of $\coneK$ is {\it exposed} if there is a nonzero $\P \in \coneK^*$ such that 
 $\coneJ = \left\{ \X \in \coneK :  \inprod{\P}{\X} = 0 \right\}$. A proper face of $\coneK$ is {\it non-exposed}, if it is not exposed. 
In general, if  $\TC(\coneJ)$ denotes the tangent linear space of 
a face $\coneJ$ of $\coneK$, then 
$\coneJ = \coneK \cap \TC(\coneJ)$. Here the {\it tangent linear space} $\TC(\coneJ)$ of 
a face $\coneJ$ of $\coneK$ is defined as the smallest linear subspace of $\spaceV$ that contains $\coneJ$. 
The {\it dimension of a face} $\coneJ$ is defined as 
the dimension of its tangent linear subspace $\TC(\coneJ)$.
We say that $\P \in \spaceV$ is {\it copositive} on a cone $\coneK \subset \spaceV$ if 
$\inprod{\P}{\X} \geq 0$ for every $\X \in \coneK$, 
{\it i.e.}, $\P \in \coneK^*$. 

Let $\H^0 \in \spaceV$. 
For every $\coneK \subset \spaceV$
 and $\rho \geq 0$, let 
$
G(\coneK,\rho)
=  \left\{ \X \in \coneK : \inprod{\H^0}{\X} = \rho\right\}. 
$
In addition, given any $\P \in \spaceV$,
we consider the following conic optimization problem 
\begin{eqnarray*}
\mbox{COP($\coneK,\P,\rho$): } \ 
\zeta(\coneK,\P,\rho)  =
 \mbox{inf} \left\{\inprod{\P}{\X} : \X \in G(\coneK,\rho)\right\} .
\end{eqnarray*}
Note that we use the convention that $\zeta(\coneK,\P,\rho)  = +\infty$ if $G(\coneK,\rho) = \emptyset$,
and that
COP($\coneK,\P,1$) coincides with COP($\coneK,\P$) introduced in Section 1. 
In the subsequent sections, we often use $\zeta(\coneK,\P,\rho)$ with $\rho \geq 0$, but COP($\coneK,\P,\rho$) 
only for $\rho=1$. 
For simplicity, we use the notation COP($\coneK,\P$) for COP($\coneK,\P,1$).

%%%%%%%%%%%%%%%%%%%%%%%%%%%%%%%%%%%
\subsection{A class of QOPs with linear equality, complementarity 
and binary constraints in nonnegative variables}

In this section,  we first consider  Burer's class   of QOPs  which were shown
to be equivalent 
to their CPP reformulations under mild assumptions (see~\eqref{eq:binary} and~\eqref{eq:complementarity} below)
in  \cite{BURER2009}.
%by Burer \cite{BURER2009}. 
 For the reader  %to g a reader, 
who might be more familiar with QOPs than POPs,
our purpose here is to show how our geometrical analysis works for QOPs, before presenting
the rigorous derivation of our convexification procedure for the POP~\eqref{eq:POP10}.

Let $\C \in \SymMat^n$, $\cc \in \Real^n$, 
$\A \in \Real^{\ell \times n}$, $\b \in \Real^{\ell}$, 
$I_{\rm bin} \subset \{1,\ldots,n\}$ (the index set for binary variables) and 
$I_{\rm comp} \subset \{(j,k): 1\leq j < k \leq n \}$ (the index set for pairs of complementary variables). 
For simplicity of notation, we assume that $I_{\rm bin} = \{1,\ldots,q\}$ for some $q \geq 0$; if $q=0$ then 
$I_{\rm bin} = \emptyset$. 
Consider a QOP of the following form:
\begin{eqnarray}
\zeta_{\mbox{\scriptsize QOP}} & = & \inf\left\{\w^T\C\w + 2\cc^T\w: 
\begin{array}{l}
\w \in \Real^n_+, \\
 f_1(\w) \equiv (\A\w - \b)^T(\A\w - \b) = \0, \\
f_2(\w)\equiv\sum_{(j,k) \in I_{\rm comp}} w_jw_k = 0, \\ 
f_{p+2}(\w) \equiv w_p(1-w_p)=0 \ (p=1,\ldots,q)
\end{array}
\right\}. 
\label{eq:QOP20}
\end{eqnarray}
Assume that the feasible region of QOP~\eqref{eq:QOP20} is nonempty.
Note that the multiple complementarity constraints $w_jw_k = 0 \ ((j,k) \in I_{\rm comp})$ in $\w \in \Real^n_+$ is written 
as the single equality constraint $f_2(\w)= 0$ in $\w \in \Real^n_+$ mainly for simplicity. 

Let 
\begin{eqnarray*}
\bGamma^{1+n} & = & \left\{\x\x^T : \x \in \Real^{1+n}_+\right\}, \ \CPP^{1+n} = \mbox{co}\bGamma^{1+n}, \\
 \CP^{1+n} & = & (\CPP^{1+n})^*  
                    = \left\{ \Y \in \SymMat^{1+n} :  \x^T\Y\x \geq 0 \ \mbox{for every } \x \in \Real^{1+n}_+ \right\}.  
\end{eqnarray*}
Then, $\bGamma^{1+n}$ forms a nonconvex cone in $\SymMat^{1+n}$. 
The convex cones $\CPP^{1+n}$ and $\CP^{1+n}$ are known as the {\it completely positive cone }
and the {\it copositive cone } in the literature \cite{BERMAN2003}, respectively. 
We know that 
\begin{eqnarray*}
\bGamma^{1+n} \subset \CPP^{1+n}  \subset \SymMat^{1+n}_+ \cap \SymN^{1+n} \subset \SymMat^{1+n}_+ \subset 
\SymMat^{1+n}_+ + \SymN^{1+n} \subset \CP^{1+n} = (\bGamma^{1+n})^*, 
\end{eqnarray*} 
where $\SymMat^{1+n}_+$ denotes the cone of positive semidefinite 
matrices in $\SymMat^{1+n}$, and $\SymN^{1+n}$ the cone of matrices with all nonnegative elements in $\SymMat^{1+n}$. 
The cone $\SymMat^{1+n}_+ \cap \SymN^{1+n}$ is often called the {\it doubly nonnegative (DNN) cone}. 

We now transform QOP~\eqref{eq:QOP20} to 
COP($\bGamma^{1+n}\cap\coneJ,\Q^0$) for some convex cone 
$\coneJ \subset \CPP^{1+n}$ and some $\Q^0 \in \SymMat^{1+n}$.
Let $m = q+2$. We first introduce 
the following homogeneous quadratic functions in $(x_0,\w) \in \Real^{1+n}$: 
\begin{eqnarray}
\left. \begin{array}{l} 
\bar{f}_0(\x) = \w^T\C\w + 2x_0\cc^T\w, \ 
 \bar{f}_1(\x) =  (\A\w - \b x_0)^T(\A\w - \b x_0), \\[5pt]
\bar{f}_2(\x) =  \sum_{(j,k) \in I_{\rm comp}}w_jw_k, \ 
\bar{f}_{p}(\x) = w_{p-2}(x_0-w_{p-2}) \ (p=3,\ldots,m). 
\end{array} \right\} \label{eq:fbar20} 
\end{eqnarray}
Then, we can rewrite QOP~\eqref{eq:QOP20} as 
\begin{eqnarray}
\zeta_{\mbox{\scriptsize QOP}} & = & \inf\left\{ \bar{f}_0(\x) : \x=(x_0,\w) \in \Real^{1+n}_+, \ x_0 = 1, \ \bar{f}_p(\x) = 0 \
(p=1,\ldots,m) \right\}. \label{eq:QOP21} 
\end{eqnarray}
Since each $\bar{f}_p(\x)$ is a homogeneous quadratic function in $\x \in \Real^{1+n}$, 
it can be rewritten as 
$\bar{f}_p(\x) = \inprod{\Q^p}{\x\x^T}$ for some $\Q^p \in \SymMat^{1+n}$ $(p=1,\ldots,n)$. 
As a result,  %we can further transform 
QOP~\eqref{eq:QOP21} can be further transformed into 
\begin{eqnarray*}
\zeta_{\mbox{\scriptsize QOP}} & = & \inf\left\{ \inprod{\Q^0}{\X} : 
\begin{array}{l}
\X \in\bGamma^{1+n}, \ \inprod{\H^0}{\X}=1, \\ 
\inprod{\Q^p}{\X} = 0 \
(p=1,\ldots,m) 
\end{array}
\right\}, 
\end{eqnarray*}
where
$\H^0$ denotes \mbox{the matrix in $\SymMat^{1+n}$ with the $(0,0)$th element $H^0_{00} = 1$}
\mbox{and 0 elsewhere.} 
We note that $\x = (1,\w) \in \Real^{1+n}_+$ iff $\x\x^T \in \bGamma$ and $\inprod{\H^0}{\X}=1$. 
By considering the convex cone 
$ 
\coneJ =\left\{ \X \in \CPP^{1+n}: \inprod{\Q^p}{\X}=0 \ (p=1,\ldots,m)\right\}, 
$ 
we can rewrite the above problem as the COP
\begin{eqnarray*}
\zeta_{\mbox{\scriptsize QOP}} & = & \inf\left\{ \inprod{\Q^0}{\X} : 
\X \in\bGamma^{1+n}\cap\coneJ, \ \inprod{\H^0}{\X}=1
\right\} \nonumber \\
& = & 
\inf\left\{ \inprod{\Q^0}{\X} : 
\X \in G(\bGamma^{1+n}\cap\coneJ,1)
\right\} 
= \zeta(\bGamma^{1+n}\cap\coneJ,\Q^0,1), 
\end{eqnarray*}
which is equivalent to QOP~\eqref{eq:QOP20}. 
Thus, we have derived COP($\bGamma^{1+n}\cap\coneJ,Q^0$) with 
a convex cone $\coneJ \subset \CPP^{1+n}$. 

If Conditions I$_0$ and II$_0$ are satisfied with $\coneK_0= \bGamma^{1+n}\cap\coneJ$, then 
COP($\coneK_0,\Q^0$) is equivalent to its covexification COP(co$\coneK_0,\Q^0$), {\it i.e.}, 
$\zeta(\coneK_0,\Q^0) = \zeta(\mbox{co}\coneK_0,\Q^0)$ (Theorem~\ref{theorem:main0}). 
If, in addition, $\coneJ$ is a face of $\CPP^{1+n}$, then co$\coneK_0 = \mbox{co}(\bGamma^{1+n}\cap\coneJ)= \coneJ$ 
(Lemma~\ref{lemma:basic}). Hence, COP($\bGamma^{1+n}\cap\coneJ,\Q^0$) is equivalent to its covexification 
COP($\coneJ,\Q^0$), which forms a CPP reformulation of QOP~\eqref{eq:QOP20} such that $\zeta(\coneJ,\Q^0,1) = \zeta_{\rm QOP}$. 

In Burer \cite{BURER2009}, the following conditions are imposed on QOP~\eqref{eq:QOP20} 
to derive its equivalent CPP reformulation:
\begin{eqnarray}
& & w_i \leq 1 \ \mbox{if } \w \in L \equiv \left\{ \w \in \Real^n_+ : \A\w - \b = \0 \right\} \ \mbox{and } i\in I_{\rm bin}, \label{eq:binary} \\
& & w_j = 0 \ \mbox{and } w_k = 0 \ \mbox{if } \w \in L_{\infty} \equiv \left\{ \w \in 
 \Real^n_+ : \A\w = \0\right\} \ \mbox{and } 
(j,k) \in I_{\rm comp}. \label{eq:complementarity}
\end{eqnarray}
%\fbox{\red{In \cite{BURER2009}, the condition \eqref{eq:complementarity} is 
%$w_j = 0$ if $\w \in L_{\infty}$ and $j\in I_{\rm bin}$. Is there a mistake here?}}\\
%\fbox{\red{\eqref{eq:complementarity} is assumed in the last paragraph of \cite[p.488]{BURER2009}. Also \eqref{eq:binary} implies the condition you}}\\
%\fbox{\red{ mentioned above.}}
Although his CPP reformulation of QOP~\eqref{eq:QOP20} is described quite differently from COP($\coneJ,\Q^0$), 
the conditions  \eqref{eq:binary} and \eqref{eq:complementarity} 
are sufficient not only for $\coneJ$ to be a face of $\CPP^{1+n}$ but also for 
$\zeta(\coneJ,\Q^0,1) = \zeta_{\rm QOP}$ to hold. This fact will be shown in Section 6.1.

%%%%%%%%%%%%%%%%%%%%%%%%%%%%%%%%%%%
\subsection{Fundamental properties of cones and their faces} 

The following lemma will play an essential role in 
the subsequent discussions. 

\lemm \label{lemma:basic20}
Let $\coneK \subset \spaceV$ be a cone. The following results hold.
\begin{description}
\item{(i) } $\coneK^* = ({\rm co}\coneK)^*$.
\item{(ii) } Assume that $\P \in \spaceV$ is copositive on $\coneK$. Then 
$\coneJ = \left\{ \X \in {\rm co}\coneK: \inprod{\P}{\X} = 0\right\}$ 
forms  an exposed face of ${\rm co}\coneK$.
\item{(iii) } Let $\coneJ_0 = {\rm co}\coneK$. Assume that 
$\coneJ_p$ is a face of $\coneJ_{p-1}$ $(p=1,\ldots,m)$. Then 
$\coneJ_{\ell}$ is a face of $\coneJ_p$ $(0 \leq p \leq \ell \leq m)$.
\end{description}
\elemm
\proof{
(i) $\coneK^* \supset ({\rm co}\coneK)^*$ follows from the fact that 
$\coneK \subset  {\rm co}\coneK$. To prove the converse inclusion, 
suppose that $\X \in \coneK^*$. Choose $\Y \in {\rm co}\coneK$ arbitrarily. 
Then there exist $\Y^i \in \coneK$ 
$(i=1,\ldots,k)$ such that 
$\Y = \sum_{i=1}^k\Y^i$. 
Since $\X \in \coneK^*$ and $\Y^i \in \coneK$, we have that 
$\inprod{\Y^i}{\X} \geq 0$ $(i=1,\ldots,k)$. 
It follows that $\inprod{\Y}{\X} = \sum_{i=1}^k\inprod{\Y^i}{\X} \geq 0$. Hence we have shown that $\inprod{\Y}{\X} \geq 0$ for 
every $\Y \in \mbox{co}\coneK$. Therefore $\X \in \mbox{co}\coneK^*$. 

(ii) Let $\X = \X^1/2+\X^2/2 \in \coneJ$, $\X^1\in{\rm co}\coneK$ and $\X^2\in{\rm co}\coneK$. By the assumption, 
$\inprod{\P}{\X^1} \geq 0$ and $\inprod{\P}{\X^2} \geq 0$. From 
$\X = \X^1/2+\X^2/2 \in \coneJ$, we also see that 
$0 = \inprod{\P}{\X} = \inprod{\P}{\X^1}/2+ \inprod{\P}{\X^2}/2$. 
Hence $\inprod{\P}{\X^1} =\inprod{\P}{\X^2} = 0$. Therefore 
$\X^1 \in \coneJ$ and $\X^2 \in \coneJ$, and we have shown that $\coneJ$ is a face 
of ${\rm co}\coneK$. Note that $\coneJ$ is exposed by definition.

(iii) We only prove the case where $m=2$ since the general case 
where $m \geq 3$ can be proved by induction. Let $\X = \X^1/2+\X^2/2 \in \coneJ_2$, $\X^1\in \coneJ_0$ and $\X^2\in\coneJ_0$. 
It follows from $\coneJ_2 \subset \coneJ_1$ that 
$\X \in \coneJ_1$. Since $\coneJ_1$ is a face of $\coneJ_0$, we obtain that $\X^1\in \coneJ_1$ and $\X^2\in \coneJ_1$. 
Now, since $\coneJ_2$ is a face of $\coneJ_1$, 
$\X^1\in \coneJ_2$ and $\X^2\in \coneJ_2$ follow. Thus we have shown that $\coneJ_2$ is a face of $\coneJ_0$. 
\qed
}

%!TEX root = ./main.tex
\section{Main results} 

Given a nonconvex cone $\coneK_0 \subset \spaceV$, 
$\H^0 \in \spaceV$ and $\Q^0 \in \spaceV$, 
the problem COP($\coneK_0,\Q^0$) minimizes the linear objective function 
$\inprod{\Q^0}{\X}$ over the nonconvex feasible region 
$G(\coneK_0,1)$. 
In Section 2.2, we have derived such a nonconvex COP from QOP~\eqref{eq:QOP20}. 
We will also see in Section 5 that a general 
class of POPs can be reformulated as such a nonconvex COP. 
By replacing $\coneK_0$  with its convex hull co$\coneK_0$, we obtain 
COP($\mbox{co}\coneK_0,\Q^0$) that minimizes the same linear 
objective function over the convex feasible region 
$G(\mbox{co}\coneK_0,1)$. Hence COP($\mbox{co}\coneK_0,\Q^0$) turns out to be a convex conic optimization problem. 
We call this process the 
{\it covexification} of COP($\coneK_0,\Q^0$). Since $\coneK_0 \subset \mbox{co}\coneK_0$, we have that 
$G(\mbox{co}\coneK_0,1) \supset G(\coneK_0,1)$ and 
$\zeta(\mbox{co}\coneK_0,\Q^0,1) \leq \zeta(\coneK_0,\Q^0,1)$ hold in general. Hence $\zeta(\mbox{co}\coneK_0,\Q^0,1)$ provides a 
lower bound for the optimal value of the original QOP or POP from which COP($\coneK_0,\Q^0$) is derived. 
If $\zeta(\mbox{co}\coneK_0,\Q^0,1) = \zeta(\coneK_0,\Q^0,1)$, we call COP($\mbox{co}\coneK_0,\Q^0$) as a {\it convex COP reformulation} of 
COP($\coneK_0,\Q^0$) (and the original QOP or POP). In this case, 
COP(co$\coneK_0,\Q^0$) solves the original QOP or POP in the sense 
that $\zeta(\mbox{co}\coneK_0,\Q^0,1)$ coincides with
their optimal values. The main result of this section is the characterization of the
convex COP reformulation of 
COP($\coneK_0,\Q^0$). 

Throughout this section, we fix a linear space $\spaceV$ 
and $\H^0 \in \spaceV$.

\subsection{A simple conic optimization problem}

For every $\coneK_0 \subset \spaceV$ and $\Q^0 \in \spaceV$,  we consider COP($\coneK_0,\Q^0$). 
To ensure $\zeta(\coneK_0,\Q^0,1) = \zeta(\mbox{co}\coneK_0,\Q^0,1)$ in Theorem~\ref{theorem:main0}, we will assume 
Conditions I$_0$ and II$_0$ introduced in Section 1. We note that 
the feasibility of COP($\coneK_0,\Q_0$) in Condition I$_0$ can be stated as $G(\coneK_0,1) \not= \emptyset$, and Condition II$_0$ 
as $\zeta(\coneK_0,\Q^0,0) \geq 0$.

Lemma~\ref{lemma:main0} and Theorem~\ref{theorem:main0} below may be regarded 
as special cases of Lemma 3.1 and Theorem 3.1 of \cite{ARIMA2018}. Although %we can derive our 
Lemma~\ref{lemma:main0} and Theorem~\ref{theorem:main0}  can be derived 
if 
$m=0$ is used in \cite{ARIMA2018}, here we present their proofs for 
the paper to be self-contained.

\lemm \label{lemma:main0} 
Let $\coneK_0 \subset \spaceV$ be a cone. Assume that Condition I$_0$ holds. Then,
\begin{description}
\item{(i) } $\mbox{\rm co}G(\coneK_0,0) = G(\mbox{\rm co}\coneK_0,0) $. 
\item{(ii) } For every $\P \in \spaceV$, $\zeta(\coneK_0,\P,0) = \zeta({\rm co}\coneK_0,\P,0)$.
\item{(iii) } For every $\P \in \spaceV$,
$ \zeta(\coneK_0,\P,0) = 
\left\{\begin{array}{ll}
0 & \mbox{if $\zeta(\coneK_0,\P,0) \geq 0$ holds,} \\
-\infty & \mbox{otherwise.}
\end{array}
\right. 
$ 
\end{description}
\elemm

\theo \label{theorem:main0} Let $\coneK_0 \subset \spaceV$ be a cone and $\Q^0 \in \spaceV$. Assume that Condition I$_0$ holds. 
Then,
\begin{description}
\item{(i) } $G(\mbox{\rm co}\coneK_0,1) =\mbox{\rm co} G(\coneK_0 ,1)+\mbox{\rm co} G(\coneK_0,0)$. 
\item{(ii) } $\zeta(\mbox{\rm co}\coneK_0,\Q^0,1) = \zeta(\coneK_0,\Q^0,1) + \zeta(\coneK_0,\Q^0,0)$.  
\item{(iii) }$\zeta(\mbox{\rm co}\coneK_0,\Q^0,1) = \zeta(\coneK_0,\Q^0,1)$  iff 
\begin{eqnarray}
	\mbox{Condition II$_0$ or \ $\zeta(\coneK_0,\Q^0,1) = -\infty$ holds.} \label{eq:iff30}
\end{eqnarray}
\end{description}
\etheo

\begin{figure}
\begin{center}
 \ifpdf
 \includegraphics[width=0.45\textwidth]{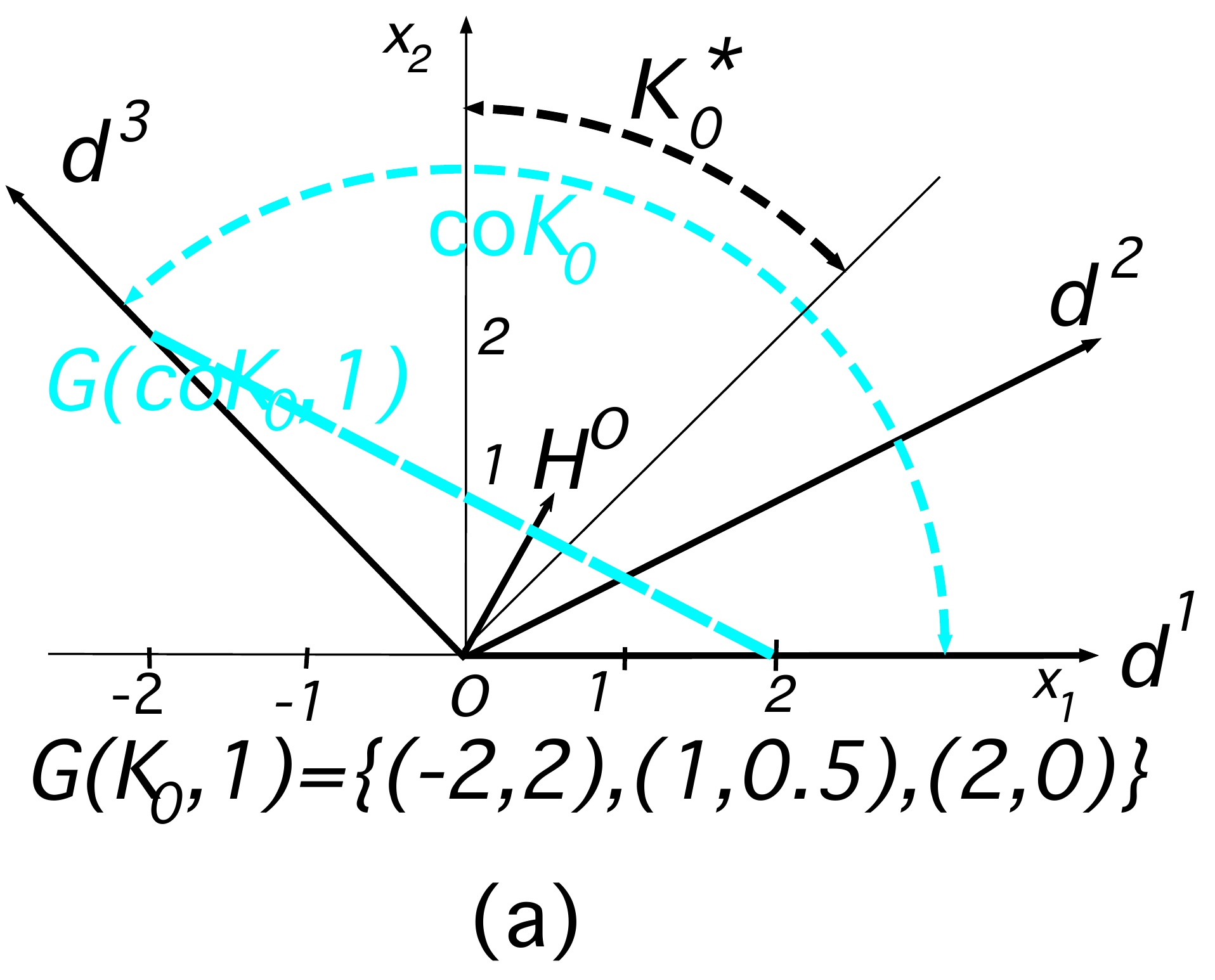}
 \else
 \includegraphics[width=0.45\textwidth]{Example1_1.eps}
 \fi
 \ifpdf
 \includegraphics[width=0.45\textwidth]{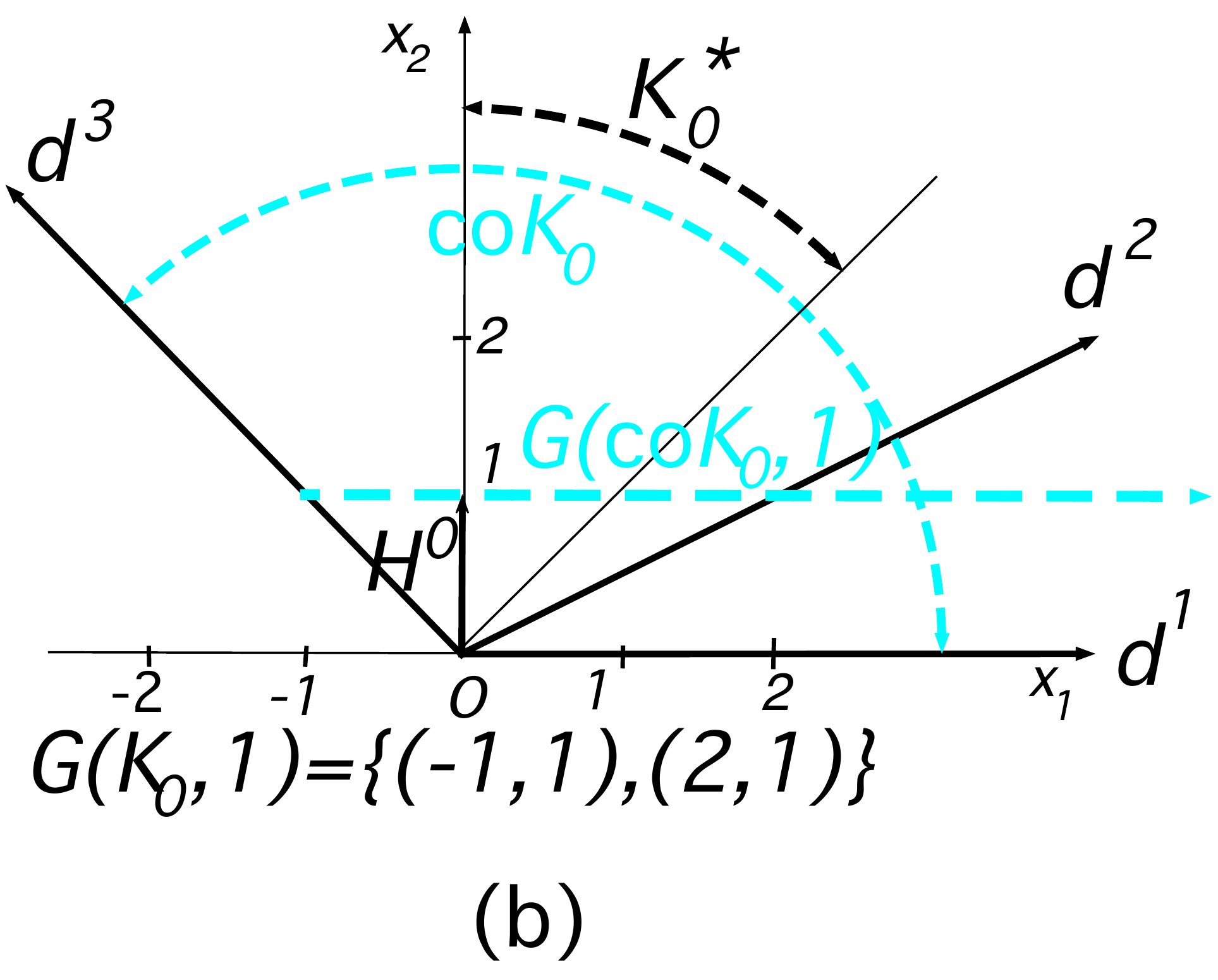}
 \else
 \includegraphics[width=0.45\textwidth]{Example1_2.eps}
 \fi
\caption{
Illustration of COP($\coneK_0,\Q^0$) and COP(co$\coneK_0,\Q^0$) under Conditions I$_0$ and II$_0$, 
where 
$\spaceV = \Real^2$, $\coneK_0 = \bigcup_{i=1}^3 \left\{\lambda\d^i: \lambda\geq 0\right\}$ and 
$G(\coneK_0,1) = \{\X \in \coneK_0 : \inprod{\H^0}{\X} = 1 \}$ (the feasible region of COP($\coneK_0,\Q^0$)). 
In case (a) where we take $\H^0 = (0.5,1) \in \coneK_0^*$, Condition I$_0$ and Condition II$_0$
are satisfied for any choice of $\Q^0 \in \Real^2$. 
In case (b) where we take $\H^0 = (0,1) \in \coneK_0^*$, Condition I$_0$ is satisfied, but Condition II$_0$ is satisfied iff 
the first coordinate $Q^0_1$ of $\Q^0 \in \Real^2$ is nonnegative. See Example~\ref{example:COP} for more details.  
}
\label{figure:Example1}
\end{center}
\end{figure}

Before presenting the proofs of Lemma~\ref{lemma:main0} and Theorem
\ref{theorem:main0}, we show an illustrative example.
%Before proving the theorem above, we present an example which illustrates the assertions of Lemma~\ref{lemma:main0} and the 
%theorem. 
\examp \label{example:COP}
Let $\spaceV=\Real^2$, $\d^1=(4,0)$, $\d^2=(4,2)$, $\d^3=(-3,3)$ and 
$
\coneK_0 = \bigcup_{i=1}^3 \left\{\lambda\d^i: \lambda\geq 0\right\}
$. 
We consider two cases (see (a) and (b) of Figure 1, respectively). 

(a) 
Let $\H^0 = (0.5,1)$, which lies in the interior of $\coneK_0^*$. In this case, we see that 
$G(\coneK_0,1) = \left\{(-2,2),(1,0.5),(2,0)\right\}$,
$G(\mbox{co}\coneK_0,1) = \mbox{co}G(\coneK_0,1) = $ the line 
segment jointing $(-2,2)$ and $(2,0)$, and $G(\coneK_0,0) = G(\mbox{co}\coneK_0,0) = \{\0\}$. 
Hence $\zeta(\mbox{co}\coneK_0,\P,0) = \zeta(\coneK_0,\P,0) = 0$ for every $\P \in \Real^2$ and Condition II$_0$ holds for every $\Q^0 \in \Real^2$. Thus 
all assertions of Lemma~\ref{lemma:main0} and Theorem~\ref{theorem:main0} hold. \vspace{-2.5mm}

\medskip
(b) Let $\H^0 = (0,1)$, which lies in the boundary of 
$\coneK_0^*$. In this case, we see that $G(\coneK_0,1) = \left\{(-1,1),(2,1)\right\}$,
$G(\mbox{co}\coneK_0,1) = \left\{(x_1,1): -1 \leq x_1 \right\}$, and $G(\coneK_0,0) = G(\mbox{co}\coneK_0,0) = \{(x_1,0): 0 \leq x_1 \}$. Hence (i) and (ii) of Lemma~\ref{lemma:main0}, and (i) 
of Theorem~\ref{theorem:main0} follow. Take 
$\Q^0 =\P=(p_1,p_2)\in \Real^2$ arbitrarily. If $p_1 \geq 0$ then 
$\zeta(\coneK_0,\P,0) = \zeta(\mbox{co}\coneK_0,\P,0) = 0$, and 
both COP($\coneK_0,\P$) and COP(co$\coneK_0,\P$) have a common 
optimal solution at $(-1,1)$ with the optimal value 
$\zeta(\coneK_0,\P,1) = \zeta(\mbox{co}\coneK_0,\P,1) = -p_1+p_2$; 
hence (iii) of Lemma~\ref{lemma:main0}, (ii) and (iii) of 
Theorem~\ref{theorem:main0} hold. Now assume that $p_1 < 0$. 
Then we see that $\zeta(\coneK_0,\Q^0,0) =\zeta(\coneK_0,\P,0) = \zeta(\mbox{co}\coneK_0,\P,0) = -\infty$. This implies that 
(ii) of Lemma~\ref{lemma:main0} holds, and that Condition II$_0$ is violated. We also see that $\zeta(\coneK_0,\Q^0,1) = 2p_1+p_2$. In this case, (iii) of Theorem~\ref{theorem:main0} asserts that $\zeta(\mbox{co}\coneK_0,\Q^0,1) \not=  \zeta(\coneK_0,\Q^0,1)$. In fact, we have that 
$-\infty = \zeta(\mbox{co}\coneK_0,\Q^0,1) < \zeta(\coneK_0,\Q^0,1) = 2p_1+p_2$. 
\eexamp

\medskip 

\noindent
{\it Proof of Lemma~\ref{lemma:main0}.} \hspace{1mm}%
(i) Since $G(\mbox{co}\coneK_0,0)$ is a convex subset of $\spaceV$ containing $G(\coneK_0,0)$, we see that 
$\mbox{co}G(\coneK_0,0) \subset G(\mbox{co}\coneK_0,0)$. 
To show the converse inclusion, assume that $\X \in G(\mbox{co}\coneK_0,0)$. Then there exist $\X^i \in \coneK_0$ 
$(i=1,2,\ldots,r)$ such that 
$\X = \sum_{i=1}^r \X^i$.
By Condition I$_0$, we know that $\inprod{\H^0}{\X^i} \geq 0$ 
$(i=1,2,\ldots,r$). Thus, each $\X^i$ % $(i=1,2,\ldots,r)$ 
satisfies 
$ \X^i \in \coneK_0 \ \mbox{and }  \inprod{\H^0}{\X^i} = 0$,  
or equivalently $\X^i \in G(\coneK_0,0)$ $(i=1,2,\ldots,r)$. Therefore, $\X = \sum_{i=1}^r \lambda_i\X^i \in \mbox{co}G(\coneK_0,0)$. 

(ii) Let $\P \in \spaceV$. We observe that 
\begin{eqnarray*}
\zeta(\coneK_0,\P,0)  
& = & \inf \left\{ \inprod{\P}{\X} : \X \in \mbox{co}G(\coneK_0,0) \right\} \ \mbox{(since $\inprod{\P}{\X}$ is linear in $\X$)} \\ 
& = & \inf \left\{ \inprod{\P}{\X} : \X \in G(\mbox{co}\coneK_0,0) \right\} \ \mbox{(by (i))} \\ 
& = & \inf \zeta(\mbox{\rm co}\coneK_0,\P,0). 
\end{eqnarray*}
 
(iii) 
Since the objective function $\inprod{\P}{\X}$ in the description of $\zeta(\coneK_0,\P,0)$
is linear and its feasible region $G(\coneK_0,0)$ 
forms a cone, we know that $\zeta(\coneK_0,\P,0) = 0$ or $-\infty$ and 
that $\zeta(\coneK_0,\P,0) = 0$ iff the objective value is nonnegative for all feasible solutions, 
{\it i.e.},  $\zeta(\coneK_0,\P,0) \geq 0$ holds. 
\qed

\medskip 

\noindent
{\it Proof of Theorem~\ref{theorem:main0}.} \hspace{1mm}
(i) 
To show the inclusion $G(\mbox{\rm co}\coneK_0,1) \subset \mbox{\rm co} G(\coneK_0,1)+\mbox{\rm co} G(\coneK_0,0)$,\ assume 
that $\X \in G(\mbox{\rm co} \coneK_0,1)$. Then there exist $\X^i \in \coneK_0 \subset \mbox{co}\coneK_0$ 
$(i=1,2,\ldots,r)$ such that 
 \begin{eqnarray*}
& &  
\X = \sum_{i=1}^r \X^i \ \mbox{and } \ 
1 =  \inprod{\H^0}{\X}  =  \sum_{i=1}^r  \inprod{\H^0}{\X^i}. 
\end{eqnarray*}
By Condition I$_0$,   $\inprod{\H^0}{\X^i} \geq 0$ $(i=1,\ldots,r)$. 
Let 
\begin{eqnarray*}
& & I_+=\left\{i: \inprod{\H^0}{\X^i}>0 \ \right\}, \ I_0=\left\{j: \inprod{\H^0}{\X^j}=0 \ \right\}, \\ 
& & \mu_i = \inprod{\H^0}{\X^i}, \ \Y^i = (1/\mu_i) \X^i  \ (i \in I_+),  \ \Y = \sum_{i \in I_+} \X^i, \\ 
& & \mu_j = 1/\left| I_0 \right|, \ \Z^j = (1/\mu_j) \X^j \  \ (j \in I_0), \ \Z = \sum_{j \in I_0} \X^j, 
\end{eqnarray*}
where $\left| I_0 \right|$ denotes the number of elements in $I_0$. Then $ \X = \Y + \Z$, and
\begin{eqnarray*}
& & \mu_i  \;> \;  0, \ \Y^i \in \coneK_0, \ 1 = \inprod{\H^0}{\Y^i} \ (i \in I_+), \ 1 
= \sum_{i \in I_+} \mu_i, \ \Y = \sum_{i\in I_+}\mu_i\Y_i , \\ 
& & \mu_j \; > \;  0, \  \Z^j \in \coneK_0, \ 0 = \inprod{\H^0}{\Z^j}  \ (j \in I_0), \  1 = \sum_{j \in I_0} \mu_j, \ \Z = \sum_{j \in I_0} \mu_j \Z^j. 
\end{eqnarray*}
Thus, 
$\Y^i \in G(\coneK_0,1)$ $(i \in I_+)$, 
$\Z^j \in G(\coneK_0,0)$ $(j \in I_0)$, $\Y \in \mbox{co}G(\coneK_0,1)$, $\Z \in \mbox{co}G(\coneK_0,0)$ and $\X = \Y+\Z$. 
Therefore, we have shown that $G((\mbox{\rm co}\coneK_0),1) \subset \mbox{co} G(\coneK_0,1)+\mbox{co} G(\coneK_0,0)$. 
In the discussion above, we have implicitly assumed that $I_0 \not= \emptyset$; otherwise $\mu_j$ $(j \in I_0)$ cannot be  consistently defined. 
If $I_0 = \emptyset$, we can just neglect $\mu_j$ and $\Z^j$ $(j \in I_0)$ and take $\Z = \O$. Then all the discussions above remain valid. 

To show the converse inclusion, suppose that $\X = \Y + \Z$ for some 
$\Y \in \mbox{co}G(\coneK_0,1)$ and $\Z \in \mbox{co} G(\coneK_0,0)$. Then we can represent 
$\Y \in \mbox{co} G(\coneK_0,1)$ as 
\begin{eqnarray*}
 \Y & = &\sum_{i=1}^p \lambda_i \Y^i, \;\; \sum_{i=1}^p \lambda_i = 1, \ \lambda_i > 0, \ \Y^i \in \coneK_0, \ \inprod{\H^0}{\Y^i} = 1 \
  (i=1,2,\ldots,p), 
\end{eqnarray*}
and $\Z \in \mbox{co} G(\coneK_0,0)$ and 
\begin{eqnarray*}
 \Z & = &\sum_{j=1}^q \lambda_i \Z^j, \ \sum_{j=1}^q \lambda_j = 1, \ \lambda_j > 0, \ \Z^j \in \coneK_0, \ \inprod{\H^0}{\Z^j} = 0 \ 
(j=1,2,\ldots,q). 
\end{eqnarray*}
Since $\mbox{co} \coneK_0$ is a  convex cone, it follows from $\Y = \sum_{i=1}^p \lambda_i \Y^i \in \mbox{co} \coneK_0$ and 
$\Z = \sum_{j=1}^q \lambda_j\Z^j \in \mbox{co} \coneK_0$ that  
$\X = \Y+\Z \in \mbox{co} \coneK_0$. We also see that 
\begin{eqnarray*}
 & & \inprod{\H^0}{\X} \;=\;
 \sum_{i =1}^p \lambda_i \inprod{\H^0}{\Y^i} + \sum_{j =1}^q \lambda_j \inprod{\H^0}{\Z^j}  =  \sum_{i =1}^p \lambda_i + 0 = 1. 
 \end{eqnarray*}
Thus,  we have shown that $\X \in G(\mbox{co}\coneK_0,1)$. 
 
(ii)
We see from (i) that 
\begin{eqnarray*}
\zeta(\mbox{co}\coneK_0,\Q^0,1) 
&= & \inf \; \left\{ \inprod{\Q^0}{\Y + \Z } : \Y \in \mbox{co}G(\coneK_0,1), \ \Z \in 
   		\mbox{co}G(\coneK_0,0)  \right\} \\
& = & \inf \; \left\{ \inprod{\Q^0}{\Y} : \Y \in \mbox{co}G(\coneK_0,1) \right\} 
+
\inf \; \left\{ \inprod{\Q^0}{\Z} : \Z \in \mbox{co}G(\coneK_0,0)  \right\} \\
& = & \inf \; \left\{ \inprod{\Q^0}{\Y} : \Y \in G(\coneK_0,1) \right\}  
+ \inf \; \left\{ \inprod{\Q^0}{\Z} : \Z \in G(\coneK_0,0)  \right\} \\
& = & \zeta(\coneK_0,\Q^0,1) + \zeta(\coneK_0,\Q^0,0).
 \end{eqnarray*}

(iii) 
``if part'': Assume that Condition II$_0$ holds. Then $\zeta(\coneK_0,\Q^0,0)= 0$ follows from Lemma~\ref{lemma:main}. 
Hence $\zeta(\mbox{co}\coneK_0,\Q^0,1) = \zeta(\coneK_0,\Q^0,1)$ by (ii). 
If $\zeta(\coneK_0,\Q^0,1)  = -\infty$, then  $\zeta(\mbox{co }\coneK_0,\Q^0,1) \leq \zeta(\coneK_0,\Q^0,1) = -\infty$. 

``only if part'': Assume that $\zeta(\mbox{co}\coneK_0,\Q^0,1) = \zeta(\coneK_0,\Q^0,1)$. 
By Condition I$_0$, $G(\coneK_0,1)$ is nonempty. Hence we have 
$\zeta(\coneK_0,\Q^0,1)  < \infty$. If $\zeta(\coneK_0,\Q^0,1)=-\infty$, then we are done. So 
suppose that $\zeta(\coneK_0,\Q^0,1)$ is finite. Then the assumption and (ii) implies that $\zeta(\coneK_0,\Q^0,0)=0$. 
Consequently, Condition II$_0$ holds.
\qed

Note that by (i) and (ii) of Lemma~\ref{lemma:main0}, we can replace $\mbox{\rm co} G(\coneK_0,0)$ and $\zeta(\coneK_0,\Q^0,0)$ in 
Theorem~\ref{theorem:main0}  by $G(\mbox{\rm co} \coneK_0 ,0)$ and $\zeta(\mbox{\rm co} \coneK_0,\Q^0,0)$, respectively. 

Next, we establish the following lemma which will play an essential role to extend %all the results 
Lemma~\ref{lemma:main0} and Theorem~\ref{theorem:main0} 
to a class of general COPs in the next section.
\lemm  \label{lemma:basic} Let $\coneK \subset \spaceV$ be a cone. Assume that $\coneJ$ is a face of {\rm co}$\coneK$. Then,
\begin{description} 
\item{(i) } {\rm co}$(\coneK\cap\coneJ) 
= \coneJ$. 
\item{(ii) } $(\coneK\cap\coneJ)^* = \coneJ^*$. 
\end{description}
\elemm
\proof{(i) 
Since $\coneJ = (\mbox{co}\coneK)\cap\coneJ$ is a convex set containing $\coneK\cap\coneJ$, we see co$(\coneK\cap\coneJ) \subset \coneJ$.
To show the converse inclusion, let $\X \in \coneJ = (\mbox{co}\coneK)\cap\coneJ$. Then there exist 
$\X^i \in \coneK \subset \mbox{co}\coneK$ 
such that 
$\X = \sum_{i=1}^m \X^i$. 
Since $\coneJ$ is a face of $\mbox{co}\coneK$, we see that $\X^i \in \coneJ$  $(i=1,\ldots,m)$. Therefore, 
$\X^i \in \coneK \cap \coneJ$  $(i=1,\ldots,m)$ and 
$\X = \sum_{i=1}^m \X^i \in \mbox{co}(\coneK\cap\coneJ)$. 

(ii) Since $(\coneK\cap\coneJ)^* =({\rm co}(\coneK\cap\coneJ))^*$ by Lemma~\ref{lemma:basic20}, 
$(\coneK \cap \coneJ)^* = \coneJ^*$ follows from (i).
\qed
}

%%%%%%%%%%%%%%%%%%%%%%%%%%%%%%%%%%%%%%
\subsection{A class of general conic optimization problems} 

For every cone 
$\coneK \subset \spaceV$, every cone $\coneJ \subset \spaceV$ and every $\Q^0 \in \spaceV$,  
we consider the class of general COPs of the following form: 
\begin{eqnarray*}
\mbox{COP($\coneK\cap\coneJ,\Q^0$): } \ 
\zeta(\coneK\cap\coneJ,\Q^0,1) &=& \inf\left\{\inprod{\Q^0}{\X} : \X \in G(\coneK\cap\coneJ,1)\right\} \nonumber\\
         &=& \inf\left\{\inprod{\Q^0}{\X} : \X \in \coneK\cap\coneJ, \inprod{\H^0}{\X} = 1\right\}. 
\end{eqnarray*}
Obviously, we can handle COP($\coneK\cap\coneJ,\Q^0$) as a special case of COP($\coneK_0,\Q^0$) by taking 
$\coneK_0 = \coneK  \cap \coneJ$. In particular, we can apply Lemma~\ref{lemma:main0} and Theorem~\ref{theorem:main0} 
if we assume Conditions I$_0$ and II$_0$ for $\coneK_0 = \coneK\cap\coneJ$. We further impose the condition that 
$\coneJ$ is a face of co$\coneK$, 
which would provide various interesting structures in COP($\coneK\cap\coneJ,\Q^0$) and a bridge between Theorem~\ref{theorem:main0} 
and many existing results on the convexification of nonconvex quadratic and polynomial optimization problems. 
By (ii) of Lemma~\ref{lemma:basic}, we know that $\coneK_0^* = (\coneK\cap\coneJ)^* = \coneJ^*$ under the assumption. 
Thus we can replace Conditions I$_0$ and II$_0$ by the following Conditions 0$_{\rm J}$, I$_{\rm J}$ and II$_{\rm J}$ for COP($\coneK\cap\coneJ,\Q^0$). 
\vspace{2mm} \\
{\bf Condition 0$_{\rm J}$:} $\coneJ$ is a face of $\mbox{co}\coneK$. \vspace{2mm} \\ 
{\bf Condition I$_{\rm J}$:} COP($\coneK\cap\coneJ,\Q^0$) is feasible, {\it i.e.}, $G(\coneK\cap\coneJ,1) \not= \emptyset$, and 
$\O \not= \H^0 \in \coneJ^* $.\vspace{0.7mm} \\ 
{\bf Condition II$_{\rm J}$:} $\inf\left\{\inprod{\Q^0}{\X} : \X \in \coneK\cap\coneJ, \ 
\inprod{\H^0}{\X}=0\right\} \geq 0$, {\it i.e.},  $\zeta(\coneK\cap\coneJ,\Q^0,0) \geq 0$. \vspace{2mm} \\ 
Note that Condition 0$_{\rm J}$ is newly added while Conditions I$_{\rm J}$ and II$_{\rm J}$ are equivalent to 
Conditions I$_0$ and II$_0$ with $\coneK_0 = \coneK \cap \coneJ$ under Condition 0$_{\rm J}$, respectively. 

Let $\coneJ$ be a face of co$\coneK$ and $\coneK_0 = \coneK\cap\coneJ$. 
Then, we know by (i) of Lemma~\ref{lemma:basic} that $\mbox{co}\coneK_0 = \mbox{co}( \coneK \cap \coneJ) = \coneJ$. 
Replacing $\coneK_0$ by $\coneK\cap\coneJ$ and $\mbox{co}\coneK_0$ by $\coneJ$ in Lemma~\ref{lemma:main0}  and in Theorem~\ref{theorem:main0}, we obtain 
the following  results in Lemma~\ref{lemma:main} and Theorem~\ref{theorem:main}.

\lemm \label{lemma:main} 
Let $\coneK \subset \spaceV$ be a cone. Assume that Conditions 0$_{\rm J}$ and I$_{\rm J}$ hold. Then,
\begin{description}
\item{(i) } $\mbox{\rm co}G(\coneK\cap\coneJ,0) = G(\coneJ,0)$. 
\item{(ii) } For every $\P \in \spaceV$, $\zeta(\coneK\cap\coneJ,\P,0) = \zeta(\coneJ,\P,0)$.
\item{(iii) } For every $\P \in \spaceV$,
$ \zeta(\coneK\cap\coneJ,\P,0) = 
\left\{\begin{array}{ll}
0 & \mbox{if Condition II$_{\rm J}$ holds}, \\
-\infty & \mbox{otherwise.}
\end{array}
\right. 
$ 
\end{description}
\elemm

\theo \label{theorem:main} Let $\coneK \subset \spaceV$ be a cone and $\Q^0 \in \spaceV$. Assume that Conditions 0$_{\rm J}$ and I$_{\rm J}$ hold. 
Then,
\begin{description}
\item{(i) } $G(\coneJ,1) =\mbox{\rm co} G(\coneK \cap\coneJ,1)+\mbox{\rm co} G(\coneK\cap\coneJ,0)$. 
\item{(ii) } $\zeta( \coneJ,\Q^0,1) = \zeta(\coneK\cap\coneJ,\Q^0,1) + \zeta(\coneK\cap\coneJ,\Q^0,0)$.  
\item{(iii) }$\zeta(\coneJ,\Q^0,1) = \zeta(\coneK\cap\coneJ,\Q^0,1)$  iff 
\begin{eqnarray}
\mbox{Condition II$_{\rm J}$ or $\zeta(\coneK\cap\coneJ,\Q^0,1) = -\infty$ holds}. \label{eq:iff31}
\end{eqnarray}
\end{description}
\etheo

The assertions of Lemma~\ref{lemma:main} and Theorem~\ref{theorem:main} are similar to but more general than 
those of Lemma 3.1 and Theorem 3.1 of \cite{ARIMA2018}, 
respectively. The essential 
difference is that our results here cover the case where $\coneJ$ can be a 
non-exposed face while those in \cite{ARIMA2018} are restricted to the case 
where $\coneJ$ is an exposed face of co$\coneK$ which is represented explicitly 
as $\coneJ = \left\{\X \in \mbox{co}\coneK : 
\inprod{\Q^p}{\X}=0 \ (p=1,\ldots,m) \right\}$ for some 
$\Q^p \in \coneK^*$ $(p=1,\ldots,m)$.

Suppose that $\coneJ$ is a face of co$\coneK$ and that its tangent space $\TC(\coneJ)$ is represented as 
$\TC(\coneJ) = \left\{ \X \in \spaceV : \inprod{\Q^p}{\X} = 0 \ (p=1,\ldots,m) \right\}$ for some 
$\Q^p \in \spaceV$ $(p=1,\ldots,m)$.  Then $\coneJ = \mbox{co}\coneK\cap\TC(\coneJ)$ and the cone $\coneJ$  is represented as in~\eqref{eq:coneJ10}. 
Therefore, COP($\coneK\cap\coneJ,\Q^0$) is equivalent to COP~\eqref{eq:COP11} introduced in Section 1. We note, however, that 
$\Q^p \in \coneK^*$ $(p=1,\ldots,m)$ may not be satisfied. 

Conversely, suppose that a COP of the form~\eqref{eq:COP11} is given. It is interesting to characterize a collection of 
$\Q^p \in \spaceV \ (p=1,\ldots,m)$ which induces a face $\coneJ$ of co$\coneK$. Such a characterization is necessary 
to construct %build up 
a class of COPs of the form~\eqref{eq:COP11} that can be reformulated as convex COPs. One sufficient condition
(which was assumed in 
\cite{ARIMA2017,ARIMA2018,KIM2013})
for $\coneJ$ defined by \eqref{eq:coneJ10}
 to be a face of co$\coneK$ is that all $\Q^p \in \spaceV \ (p=1,\ldots,m)$ 
are copositive on $\coneK$. However, this sufficient condition can sometimes be restrictive. For example, we can replace $\Q^m$ by 
$-\sum_{p=1}^m\Q^p$ to generate the same $\coneJ$ 
but $-\sum_{p=1}^m\Q^p$ is no longer copositive on $\coneK$.  We also see that this sufficient condition ensures that $\coneJ$ defined 
by~\eqref{eq:coneJ10} is an exposed face of co$\coneK$. In fact, in this case, $\coneJ$ coincides with 
$\left\{ \X \in \mbox{co}\coneK : \sum_{p=1}^m\inprod{\Q^p}{\X} = 0  \right\}$. 
If $\coneJ$ is a non-exposed face of $\mbox{co}\coneK$, 
%we cannot represent
 $\coneJ$ cannot be represented in terms of any collection of copositive $\Q^p$ on $\coneK$ $(p=1,\ldots,m)$ as in~\eqref{eq:coneJ10}. 
We will investigate such cases in Section~4.

%!TEX root = ./main.tex
\section{Copositivity conditions} 

Throughout this section, we fix a linear space $\spaceV$, 
a cone $\coneK \subset \spaceV$ and $\H^0 \in \spaceV$. 
In Section 3.2, we have shown that if $\coneJ$ is a face of co$\coneK$, we can always represent $\coneJ$ as in~\eqref{eq:coneJ10}
for some $\Q^p \in \spaceV$ $(p=1,\ldots,m)$. 
In Section 4.1, we strengthen this equivalence relation by introducing a  {\it hierarchy of copositivity condition} and show how 
we can choose such $\Q^p \in \spaceV$ $(p=1,\ldots,m)$ to satisfy the condition 
recursively. 
The hierarchy of copositivity condition was originally 
proposed in Arima, Kim and Kojima \cite{ARIMA2012} as a condition for characterizing a class of QOPs that are equivalent to their CPP
reformulations. Here, we extend the condition to a more general class of COP($\coneK\cap\coneJ,\Q^0$), which includes their class of QOPs.
In Section 4.2, we present some characterizations of the 
copositivity of $\P \in \spaceV$ on a face $\coneJ$ of co$\coneK$.  They are useful to construct a face $\coneJ$ of co$\coneK$ 
in terms of $\Q^p \in \spaceV$ $(p=1,\ldots,m)$ as in~\eqref{eq:coneJ10}. 

\subsection{The hierarchy of copositivity condition}

Recall that $\coneJ$ is an exposed face of co$\coneK$ iff $\coneJ = \left\{ \X \in \mbox{co}\coneK : 
\inprod{\Q^1}{\X} = 0 \right\}$ for some copositive $\Q^1 \in  \spaceV$ on $\mbox{co}\coneK$, 
{\it i.e.}, $\Q^1 \in (\mbox{co}\coneK)^*$. 
The two  lemmas  below generalize this fact, assuming implicitly that 
$\coneJ$ can be a non-exposed face of $\mbox{co}\coneK$.  
(As we have mentioned in Section 1, some faces of 
the CPP cone $\CPP^{1+n}$ are non-exposed if $n \geq 5$ \cite{ZHANG2018}.)

\lemm \label{lemma:face40}
Let $\coneJ$ be a proper face of $\mbox{co}\coneK$. Let $\coneJ_0=\mbox{co}\coneK$. Then there exist sequences of faces 
$\coneJ_1,\ldots,\coneJ_m$ of $\mbox{co}\coneK$ and 
$\Q^1,\ldots,\Q^m \in \spaceV$ for some positive integer $m$ such that 
\begin{eqnarray}
\left. 
\begin{array}{l}
\O\not=\Q^p \in \coneJ_{p-1}^*\cap\TC(\coneJ_{p-1}), \\[5pt]
 \coneJ_p  = \left\{ \X \in \coneJ_{p-1} : 
\inprod{\Q^p}{\X} = 0 \right\}, 
 \mbox{\rm dim}\coneJ_{p-1} > \mbox{\rm dim}\coneJ_p 
\ \mbox{and } \coneJ_m = \coneJ 
\end{array}
\right\} 
\label{eq:hierarchy40}
\end{eqnarray}
$(p=1,\ldots,m)$. 
\elemm
\proof{Let $\overline{\X}$ be a relative interior point of $\coneJ$ with 
respect to the tangent space $\TC(\coneJ)$ of $\coneJ$. Since $\overline{\X}$ is a boundary point 
of the cone $\coneJ_0 = \mbox{co}\coneK$ with respect to the tangent space $\TC(\coneJ_0)$, 
we can take a supporting hyperplane of $\coneJ_0$ at $\overline{\X}$ in the tangent space $\TC(\coneJ_0)$, say, 
$\left\{ \X \in \TC(\coneJ_0) : \inprod{\Q^1}{\X} = 0\right\}$ for 
some nonzero $\Q^1 \in \coneJ_0^{*}\cap\TC(\coneJ_0)$. Let $\coneJ_1 = \left\{ \X \in \coneJ_0:  \inprod{\Q^1}{\X} = 0 \right\}$, 
which forms a face of $\coneJ_0 = \mbox{co}\coneK$ by (ii) of Lemma~\ref{lemma:basic20}. % , hence, a face of $\mbox{co}\coneK$ by (iii) of 
% Lemma~\ref{lemma:basic20}.  
By construction, 
$\coneJ \subset \coneJ_1 \subset \coneJ_0$ and $\mbox{dim}\coneJ \leq \mbox{dim}\coneJ_1 < \mbox{dim}\coneJ_0$. If $\coneJ = \coneJ_1$, 
$\TC(\coneJ) = \TC(\coneJ_1)$ or  $\overline{\X}$ lies in the relative interior of $\coneJ_1$ with respect to $\TC(\coneJ_1)$, 
we are done. In general, suppose that $\overline{\X}$ is a 
relative boundary point of a face $\coneJ_{p-1}$ with 
respect to $\TC(\coneJ_{p-1})$ $(1 \leq p)$, we can take a supporting hyperplane of $\coneJ_{p-1}$ at 
$\overline{\X}$ in the tangent space $\TC(\coneJ_{p-1})$, say $\left\{ \X \in \TC(\coneJ_{p-1}) : \inprod{\Q^p}{\X} = 0\right\}$ for 
some nonzero $\Q^p \in \coneJ_{p-1}^{*}\cap\TC(\coneJ_{p-1})$. Let $\coneJ_p = \left\{ \X \in \coneJ_{p-1} : \inprod{\Q^p}{\X} = 0\right\}$. 
Since $\mbox{dim}\coneJ_{p-1} > \mbox{dim} \coneJ_p$, this process terminates in a finite number of steps to obtain 
 a sequence of faces 
$\coneJ_1,\ldots,\coneJ_m$ of $\mbox{co}\coneK$ and a sequence 
$\Q^1,\ldots,\Q^m \in \spaceV$ satisfying \eqref{eq:hierarchy40}.
\qed
}

 Note that Lemma   \ref{lemma:face40}  shows that any proper face $\coneJ$ of co$\coneK$ can be represented in terms of 
{\it a hierarchy of copositivity condition}: 
\begin{eqnarray}
\coneJ_0 & = & \mbox{co}\coneK,  \label{eq:hierarchy41} \\ 
\coneJ_p & = & \left\{ \X \in \coneJ_{p-1} : \inprod{\Q^{pj}}{\X} = 0 \ (j=1,\ldots,q_p) \right\} \nonumber \\
& \ &  \mbox{for some copositive $\Q^{pj} \in \spaceV  \ (j=1,\ldots,q_p)$  on $\coneJ_{p-1}$} 
                 \ (p=1,\ldots,m), \label{eq:hierarchy42} \\ 
 \coneJ & = & \coneJ_m  =  \left\{ \X \in \mbox{co}\coneK : \inprod{\Q^{pj}}{\X} = 0 \ (j=1,\ldots,q_p,p=1,\ldots,m) \right\}
\label{eq:hierarchy43}
\end{eqnarray}
for some positive integers $q_p$ $(p=1,\ldots,m)$ and  $m$. 

Conversely, we can construct 
any face $\coneJ$ of co$\coneK$ by \eqref{eq:hierarchy41}, \eqref{eq:hierarchy42} and \eqref{eq:hierarchy43} as 
we shall present next.
Since all $\Q^{pj} \in \spaceV  \ (j=1,\ldots,q_p)$ 
are copositive on $\coneJ_{p-1}$ in \eqref{eq:hierarchy42}, we can replace \eqref{eq:hierarchy42} by 
\begin{eqnarray}
\coneJ_p & = & \left\{ \X \in \coneJ_{p-1} : \inprod{\Q^{p}}{\X} = 0 \right\} \nonumber \\
& \ &  \mbox{for some copositive $\Q^{p} \in \spaceV $  on $\coneJ_{p-1}$} \ (p=1,\ldots,m)  \label{eq:hierarchy44}
\end{eqnarray}
as in Lemma~\ref{lemma:face40} by letting $\Q^p = \sum_{j=1}^{q_p} \Q^{qj}$. 
We also see that if $0 \leq k < p \leq m$ and $\P \in \spaceV$ is copositive on $\coneJ_k$, then it is copositive on $\coneJ_p$ since 
$(\coneJ_k)^* \supset (\coneJ_p)^*$. 
This implies that replacing 
\eqref{eq:hierarchy42} by \eqref{eq:hierarchy44} is not restrictive at all. 
Furthermore, if $\coneJ_{p-1}$ is a face of 
co$\coneK$, then $\coneJ_{p-1} = \mbox{co}(\coneK\cap\coneJ_{p-1})$ 
by (i) of Lemma~\ref{lemma:basic}. Hence, % we can replace 
``copositive on $\coneJ_{p-1}$'' can be replaced by ``copositive on 
$\coneK\cap\coneJ_{p-1}$'' in \eqref{eq:hierarchy42} and \eqref{eq:hierarchy44}.

\lemm \label{lemma:face41}
Let $\coneK \subset \spaceV$ be a cone. 
Let $\Q^p \in\spaceV$ $(p=0,\ldots,m)$ be given, and construct a sequence of $\coneJ_p \subset \spaceV$ $(p=0,\ldots,m)$ by 
\begin{eqnarray}
\coneJ_0 = \mbox{co}\coneK \ \mbox{and } \coneJ_p = \left\{ \X \in \coneJ_{p-1} : \inprod{\Q^{p}}{\X} = 0 \right\} \ (p=1,\ldots,m) 
\label{eq:Jp40}
\end{eqnarray}
Assume that $\Q^{p} \in \spaceV$ is copositive on 
$\coneK\cap\coneJ_{p-1}$ $(p=1,\ldots,m)$. Then each $\coneJ_{p}$ 
is a face of $\coneJ_{p-1}$ and a face of co$\coneK$ $(p=1,\ldots,m)$. 
\elemm 
\proof{The assertion follows from (ii) and (iii) of Lemma~\ref{lemma:basic20}. 
%We only prove the case where $m=2$ since the general case 
%where $m \geq 3$ can be proved by induction. 
%By (ii) of Lemma~\ref{lemma:basic20}, we know that $\coneJ_1$ is a face of 
%$\coneJ_{0}$ and that $\coneJ_2$ is a face of 
%$\coneJ_{1}$.  
%Let $\X = \X^1/2+\X^2/2 \in \coneJ_2$, $\X^1\in {\rm co}\coneK$ and $\X^2\in{\rm co}\coneK$. It follows from $\coneJ_2 \subset \coneJ_1$ that 
%$\X \in \coneJ_1$. Since $\coneJ_1$ is a face of ${\rm co}\coneK$, we obtain that $\X^1\in \coneJ_1$ and $\X^2\in \coneJ_1$. 
%Now, since $\coneJ_2$ is a face of $\coneJ_1$, 
%$\X^1\in \coneJ_2$ and $\X^2\in \coneJ_2$ follow.
%By (iv) of Lemma~\ref{lemma:basic20}, $\coneJ=\coneJ_m$ is a 
%a face of $\coneK$. 
\qed
} 

It should be noted that 
$\Q^{1}$ need to be chosen from the cone 
$\coneK^*$, but $\Q^{p}$ from a possibly wider cone $\coneJ_{p-1}^*$ 
than $\coneJ_{p-2}^*$ $(p=2,\ldots,m)$. 

\subsection{Characterization of copositivity}

Let $\coneJ_0 = \mbox{co}\coneK$.  
We assume that $k=0$ or a face $\coneJ_k$ of $\coneJ_{k-1}$ has already  been constructed through~\eqref{eq:Jp40} 
for some $k=1,\ldots,p-1$. Now, we focus %our attention 
on the choice of a copositive 
$\Q^{p}   \in \spaceV$ on $\coneK\cap\coneJ_{p-1}$, {\it i.e.}, $\Q^p \in (\coneK\cap\coneJ_{p-1})^*$ so that the cone 
$
\coneJ_{p} = \left\{ \X \in \coneJ_{p-1} : \inprod{\Q^{p}}{\X}=0 \right\}
$
can become a face of $\coneJ_{p-1}$. 
By definition, $\Q^p$ is copositive on $\coneK\cap\coneJ_{p-1}$ iff 
\begin{eqnarray}
&& \eta(\coneK\cap\coneJ_{p-1},\Q^p) \equiv  \inf\left\{ \inprod{\Q^p}{\X} : \X \in \coneK\cap\coneJ_{p-1} \right\} \geq 0.
\label{eq:copositive0} 
\end{eqnarray}

\lemm \label{lemma:copositivity2} Let $\H^0 \in \coneK^*$ and 
$\coneJ_{p-1}$ be a face of co$\coneK$. 
Assume that $\zeta(\coneK\cap\coneJ_{p-1},\Q^p,0) \geq 0$.  Then, 
\eqref{eq:copositive0} is equivalent to either of the following 
two conditions:
\begin{eqnarray}
&& \zeta (\coneK\cap\coneJ_{p-1},\Q^p,\rho) \geq 0 \ \mbox{for every $\rho \geq 0$}, \label{eq:copositive2} \\[5pt]
&& \zeta(\coneK\cap\coneJ_{p-1},\Q^p,1) \geq 0.  \label{eq:copositive3}
\end{eqnarray}
\elemm
\proof{
Since 
$ \coneK\cap\coneJ_{p-1} \supset G(\coneK\cap\coneJ_{p-1},\rho) \ \mbox{for every $\rho \geq 0$}$, 
we see that \eqref{eq:copositive0} $\Rightarrow$ 
\eqref{eq:copositive2} $\Rightarrow$ \eqref{eq:copositive3}.
Thus, it suffices to show that 
\eqref{eq:copositive3} $\Rightarrow$ \eqref{eq:copositive2} 
$\Rightarrow$ \eqref{eq:copositive0}. 
 
\eqref{eq:copositive3} $\Rightarrow$ \eqref{eq:copositive2}:  
Assume that \eqref{eq:copositive3} holds.
Let $\X \in \coneK\cap\coneJ_{p-1}$ and $\inprod{\H^0}{\X} = \rho$. 
First, we consider the case $\rho > 0$. Then, $\inprod{\H^0}{\X/\rho} = 1$ and $\X/\rho \in  \coneK\cap\coneJ_{p-1}$. As a result,
$\inprod{\Q^p}{\X/\rho} \geq 0$, which implies that $\inprod{\Q^p}{\X} \geq 0$. Therefore $\zeta(\coneK\cap\coneJ_{p-1},\Q^p,\rho) \geq 0$. 
The second case where $\rho=0$ simply  follows from the assumption that $\zeta(\coneK\cap\coneJ_{p-1},\Q^p,0) \geq 0$.  

\eqref{eq:copositive2} $\Rightarrow$ \eqref{eq:copositive0}: Assume that \eqref{eq:copositive2} holds. Take $\X \in \coneK\cap\coneJ_{p-1}$ arbitrarily.  
It follows from $\X \in \coneJ_{p-1}$ and $\H^0  \in \coneK^* \subset \coneJ_{p-1}^*$ that $\rho = \inprod{\H^0}{\X} \geq 0$. Hence 
$\X \in G(\coneK\cap\coneJ_{p-1},\rho)$ with $\rho \geq 0$. Thus $\inprod{\Q^p}{\X} \geq 0$ follows from~\eqref{eq:copositive2}.  
\qed
} 

\rema 
Suppose that \eqref{eq:copositive3} holds. 
If $G(\coneK\cap\coneJ_{p},1) \not= \emptyset$, then $\inprod{\Q^p}{\X} = 0$ for some $\X \in G(\coneK\cap\coneJ_{p-1},1)$. Thus,
$\zeta(\coneK\cap\coneJ_{p-1},\Q^p,1) = 0$ and  $\zeta(\coneK\cap\coneJ_{p-1},\Q^p,\rho) = 0$ for every $\rho > 0$. 
This implies that $\lim_{\rho \rightarrow 0+} \zeta(\coneK\cap\coneJ_{p-1},\Q^p,\rho) = 0$. 
By (iii) of Lemma~\ref{lemma:main}, we also know 
that either $\zeta(\coneK\cap\coneJ_{p-1},\Q^p,0) = 0$ or $\zeta(\coneK\cap\coneJ_{p-1},\Q^p,0) = -\infty$. 
Thus the assumption made in  Lemma \ref{lemma:copositivity2}  is to ensure that the latte case where 
$\zeta(\coneK\cap\coneJ_{p-1},\Q^p,0) = -\infty < \lim_{\rho \rightarrow 0+} \zeta(\coneK\cap\coneJ_{p-1},\Q^p,\rho) = 0$ 
(a discontinuity of  $\zeta(\coneK\cap\coneJ_{p-1},\Q^p,\rho)$ at $\rho=0$) cannot occur. 
Note that  if $G(\coneK\cap\coneJ_{p-1},0) = \{\0\}$, then clearly 
$\zeta(\coneK\cap\coneJ_{p-1},\Q^p,0) = 0$, and the assumption of the lemma holds.
\erema

\rema 
We consider the case where $\coneJ_{p-1}$ is an exposed face of $\mbox{co}\coneK$, where there is a nonzero $\H^1 \in \coneK^{*}$ 
such that $\coneJ_{p-1} = \mbox{co}\coneK \cap \spaceL$ with $\spaceL = \left\{ \X \in \spaceV : \inprod{\H^1}{\X} = 0 \right\}$. 
Then $\coneJ_{p-1}^* = \mbox{cl}(\coneK^* + \spaceL^{\perp})$. Now assume that $\Q^p \in \spaceV$ is copositive 
on $\coneK\cap\coneJ_{p-1}$ or 
$\Q^p \in (\coneK\cap\coneJ_{p-1})^* = \coneJ_{p-1}^*$. If $\coneK^* + \spaceL^{\perp}$ is closed, then $\coneJ_{p-1}^* = \coneK^* + \spaceL^{\perp}$. Hence there exist $\widehat{\Y} \in \coneK^*$ 
and $\hat{y}_1 \in \Real$ such that $\Q^p = \widehat{\Y} + \H^1\hat{y_1}$. It follows that 
\begin{eqnarray*}
\coneJ_{p} & = & \left\{ \X \in \mbox{co}\coneK :  \inprod{\H^1}{\X} = 0, \ \inprod{\Q^p}{\X} = 0 \right\}  \\ 
            & = &  \left\{ \X \in \mbox{co}\coneK :  \inprod{\H^1}{\X} = 0, \ \inprod{\widehat{\Y}+\H^1\hat{y}_1}{\X} = 0 \right\}  \\ 
            & = &  \left\{ \X \in \mbox{co}\coneK :  \inprod{\H^1}{\X} = 0, \ \inprod{\widehat{\Y}}{\X} = 0 \right\} \\ 
            & = & \left\{ \X \in \mbox{co}\coneK :  \inprod{\H^1+\widehat{\Y}}{\X} = 0 \right\}. \ \mbox{(since $\H^1, \ \widehat{\Y} \in \coneK^{*}$)} 
\end{eqnarray*}
This implies that $\coneJ_{p}$ is also an exposed face of co$\coneK$. 
\erema

%!TEX root = ./main.tex
%\red{
\section{Convex COP reformulation of polynomial optimization problems}
%}

%\red{

%The purpose of this section is to
We extend the CPP reformulation of 
QOPs studied in many papers such as \cite{BOMZE2002,BOMZE2017,BURER2009} (see also Sections 2.2 and \red{6.1}.) 
to POPs.
The results presented in this section are closely related to those in Section~3 of 
\cite{ARIMA2018b}, but 
our class of POPs of the form~\eqref{eq:POP10} that can be reformulated as convex COPs does cover  POPs in nonnegative variables with 
polynomial equality constraints satisfying the hierarchy of copositivity conditions, which is more general than 
the copositivity condition assumed in Section 3 of \cite{ARIMA2018b}. 

To apply the results  described in Sections~3 and~4 to a convex conic reformulation of POP~\eqref{eq:POP10}, 
we first reduce POP~\eqref{eq:POP10} to COP($\bGamma^{\sAC}\cap\coneJ,\Q^0$). 
Here $\bGamma^{\sAC}$ is a nonconvex cone in a linear space $\SymMat^{\sAC}$ of symmetric matrices  
whose dimension depends on the maximum degree of
the monomials involved in $f_i(\w)$ $(i=0,\ldots,m)$ of POP~\eqref{eq:POP10}, and
$\AC$ stands for a set of monomials. 
The convex hull of $\bGamma^{\sAC}$, denoted as $\CPP^{\sAC}$,  corresponds to an extension of  
the CPP cone $\CPP^{1+n}$. The polynomial function $f_p(\w)$ is converted into $\inprod{\Q^p}{\X}$ in $\X \in \bGamma^{\sAC}$ 
for some $\Q^p \in \SymMat^{\sAC}$
with the additional constraint $\inprod{\H^0}{\X} = 1$ 
through its homogenization $\bar{f}_p(\x)$  $(p=0,\ldots,m)$, 
and then the face $\coneJ$ of $\CPP^{\sAC}$ is defined as in \eqref{eq:coneJ10}. 

We explain how  a polynomial function in $\w \in \Real^n$ is homogenized in Section 5.1, and 
define an extended completely positive cone $\CPP^{\sAC}$ in Section 5.2. The conversion of POP~\eqref{eq:POP10} into 
COP($\bGamma^{\sAC}\cap\coneJ,\Q^0)$ is presented in Section 5.3, 
and the convex reformulation of COP($\bGamma^{\sAC}\cap\coneJ,\Q^0)$ into COP($\coneJ,\Q^0)$ is discussed in Section 5.4. 

%%%%%%%%%%%%%%%%%%%%%%%%%%%%%
\subsection{Homogenizing polynomial functions} 

Let $\tau$ be a positive integer. 
We call that a real valued polynomial function $\bar{f}(\x)$ in $\x \in \Real^{1+n}$ is {\it homogeneous } with degree $\tau \in \Integer_+$ (or degree $\tau$  homogeneous) 
if $\bar{f}(\lambda\x) = \lambda^\tau \bar{f}(\x)$ for every $\x \in \Real^{1+n}$ and $\lambda \geq 0$. 
For the consistency of the discussions throughout Section~5, a homogeneous polynomial function is 
defined in $\Real^{1+n}$ but not $\Real^n$, where the first coordinate of $\Real^{1+n}$ is
indexed by $0$;  
we write $\x = (x_0,x_1,\ldots,x_n)$ or $\x = (x_0,\w)$ with $\w \in \Real^n$. 

For each $\balpha = (\alpha_1,\ldots,\alpha_n) \in \Integer^{n}_+$, let $\w^{\salpha}$ denote the monomial 
$\prod_{i=1}^nw_i^{\alpha_i}$ with degree $\tau_0 = \left| \balpha \right| \equiv \sum_{i=1}^n \alpha_i$. 
Let $\tau$ be a nonnegative integer no less than $\tau_0$. By introducing an 
additional variable $x_0 \in \Real$, which will be fixed to $1$ later, 
we can convert 
the previous monomial to the monomial $x_0^{\tau-\tau_0}\w^{\salpha}$ in 
$(x_0,w_1,\ldots,w_n) \in \Real^{1+n}$ with degree $\tau$. Using this technique, 
we can convert any polynomial function $f(\w)$ in $\w = (w_1,\ldots,w_n) \in \Real^n$ with degree $\tau_0$ to a homogeneous 
polynomial function $\bar{f}(x_0,\w)$ in $(x_0,w_1,\ldots,w_n)  \in \Real^{1+n}$ with degree $\tau \geq \tau_0$ such that $\bar{f}(1,\w) =f(\w)$ for 
every $\w = (w_1,\ldots,w_n) \in \Real^n$. 

%%%%%%%%%%%%%%%%%%%%%%%%%%%%%%%%
\subsection{An extension of the completely positive cone}

We begin by introducing some additional notation and symbols.
For each positive integer $\omega$, 
we define $\AC_{\omega} = \left\{\balpha=(\alpha_0,\ldots,\alpha_n) \in \Integer^{1+n}_+ :  \left| \balpha \right| = \omega \right\}$.
For each nonempty subset $\AC$ of $\AC_{\omega}$, let
$\Real^{\sAC}$ be the $\left| \AC \right|$-dimensional Euclidean space 
whose coordinates are indexed by 
$\balpha \in \AC$, where $\left| \AC \right|$ stands for the cardinality of 
$\AC$, {\it i.e.}, the number of elements in $\AC$.   
We use $\SymMat^{\sAC} \subset \Real^{\sAC} \times  \Real^{\sAC} $ to denote the space of 
$\left| \AC \right| \times \left| \AC \right|$ symmetric matrices whose elements are indexed by 
$(\balpha,\bbeta) \in  \AC  \times  \AC$. 
Let $\SymMat^{\sAC}_+$ denote the cone of positive semidefinite matrices in $\SymMat^{\sAC}$, and 
$\SymN^{\sAC}$ the cone of nonnegative matrices  in $\SymMat^{\sAC}$. 

Let $\omega$ be a  positive integer and $\emptyset \not= \AC \subset \AC_{\omega}$. 
We define 
\begin{eqnarray*}
\bGamma^{\sAC} &=& \left\{ \u^{\sAC}(\x)(\u^{\sAC}(\x))^T \in \SymMat^{\sAC}: \x \in \Real^{1+n}_+\right\} \ \mbox{and} \ 
\CPP^{\sAC} = \mbox{co}\bGamma^{\sAC}. 
\end{eqnarray*}
Here $\u^{\sAC}(\x)$ denotes the $\left|\AC\right|$-dimensional column vector 
of monomials $\x^{\salpha}$ $(\balpha \in \AC)$.
We note that every element $[\u^{\sAC}(\x)(\u^{\sAC}(\x))^T] _{\salpha\sbeta}= \x^{\salpha}\x^{\sbeta}$ 
is a degree $2\omega$ monomial in $\x \in \Real^{1+n}$. It follows that $\bGamma^{\sAC}$ forms a cone in $\SymMat^{\sAC}$.
The coordinate indices $(\balpha \in \AC)$ 
are ordered so that $\u^{\sAC}(\x) \in \Real^{\sAC}$ for every $\x \in \Real^{1+n}$.
We call $\CPP^{\sAC}$ an
{\it extended completely positive cone}, and the dual of  $\CPP^{\sAC}$, 
$\CP^{\sAC}  = (\CPP^{\sAC})^* = (\bGamma^{\sAC})^*$ an {\it extended copositive cone}. 

By definition, we know that $\bGamma^{\sAC} \subset  \SymMat^{\sAC} \cap \SymN^{\sAC}$. 
We also observe that 
\begin{eqnarray*}
X_{\salpha\sbeta} = \x^{\salpha + \sbeta} = \x^{\sgamma+\sdelta} = 
X_{\sgamma\sdelta} \ \mbox{if } \balpha + \bbeta = \bgamma+\bdelta
\end{eqnarray*}
for every $\X = \u^{\sAC}(\x)(\u^{\sAC}(\x))^T \in \bGamma^{\sAC}$. 
This implies that the cone $\bGamma^{\sAC}$ and its convex hull $\CPP^{\sAC}$ is contained in the linear subspace 
$\spaceL^{\sAC}$ of $\SymMat^{\sAC}$ defined by 
\begin{eqnarray*}
\spaceL^{\sAC} & = & \left\{ \X \in \SymMat^{\sAC} :  X_{\salpha\sbeta} = X_{\sgamma\sdelta} \ \mbox{if } \balpha + \bbeta = \bgamma+\bdelta \right\}. 
\end{eqnarray*}
Therefore, 
\begin{eqnarray}
 \bGamma^{\sAC} \subset 
\CPP^{\sAC} \subset \SymMat^{\sAC}_+\cap\SymN^{\sAC}\cap\spaceL^{\sAC}
\subset \SymMat^{\sAC}_+ \subset  \SymMat^{\sAC}_+ + \SymN^{\sAC} 
+ \left(\spaceL^{\sAC}\right)^{\perp} \subset \CP^{\sAC} = (\bGamma^{\sAC})^*.
\label{eq:inclusion50}
\end{eqnarray}

Let $\bar{f}(\x)$ be a degree $2\omega$ homogeneous polynomial 
function. Then, we can write $\bar{f}(\x) = \sum_{\sgamma \in \sBC}c_{\sgamma}\x^{\sgamma}$ for 
some nonzero $c_{\sgamma} \in \Real$ $(\bgamma \in \BC)$ and some $\BC \subset \AC_{2\omega}$. 
Since 
$\AC_{\omega} + \AC_{\omega} \equiv \{\balpha+\bbeta : \balpha \in \AC_{\omega}, \ \bbeta \in  \AC_{\omega} \} 
= \AC_{2\omega} \supset \BC$
and the matrix $\u^{\sAC}(\x)(\u^{\sAC}(\x))^T \in \SymMat^{\sAC}$ involves all monomials in $\AC + \AC$ 
for every $\AC \subset \AC_{\omega}$,  
we can choose an $\AC \subset \AC_{\omega}$ such that $\BC \subset \AC + \AC$ (see \cite{KOJIMA2005} for such a 
choice $\AC$ from $\AC_{\omega}$),  and  a matrix $\P \in \SymMat^{\sAC}$ such that 
$\bar{f}(\x) = \inprod{\P}{\u^{\sAC}(\x)(\u^{\sAC}(\x))^T}$ for every $\x \in \Real^{1+n}$. 
(Note that such a $\P \in \SymMat^{\sAC}$ is not unique.)
In our subsequent discussion, we impose an additional condition that $\AC$ contains $\balpha^{\omega} \equiv (\omega,0,\ldots,0) \in \Real^{n}$, 
and assume that the first coordinate of $\Real^{\sAC}$ is $\balpha^{\omega}$, the upper-leftmost element of each $\X \in \SymMat^{\sAC} \subset 
\Real^{\sAC}\times\Real^{\sAC}$ is $X_{\salpha^{\omega}\salpha^{\omega}}$ and that the first element of $\u^{\sAC}(\x) \in \Real^{\sAC}$ is $\x^{\salpha^{\omega}} = 
x_0^{\omega}$. 

As a consequence of the representation of $\bar{f}(\x) = \inprod{\P}{\u^{\sAC}(\x)\u^{\sAC}(\x)^T}$, 
$\bar{f}(\x) \geq 0$ for every 
$\x \in \Real^{1+n}_+$ iff $\inprod{\P}{\X} \geq 0$ for every $\X \in \bGamma^{\sAC}$ or equivalently 
$\P \in (\bGamma^{\sAC})^* = \CP^{\sAC}$.  

\subsection{Conversion of POP~\eqref{eq:POP10} to COP$(\bGamma^{\sAC}\cap\coneJ,\Q^0)$} 

Let 
$\tau_{\min} = \max\{\mbox{deg}f_i(\w) : i=0,\ldots,m\}$ and $\omega$ be a positive integer such that 
$2\omega \geq \tau_{\min}$. 
By applying the homogenization technique with degree $2\omega$ described in the previous section to the polynomial 
function $f_i(\x)$ $(i=0,\ldots,m)$, we can convert POP~\eqref{eq:POP10} to 
\begin{eqnarray}
\zeta^* =\inf \left\{ \bar{f}_0(\x) : \x = (x_0,\w) \in \Real^{1+n}, \ \bar{f}_i(\x)=0 \ (i=1,\ldots,m), \ x_0 = 1 \right\}.  
\label{eq:POP51}
\end{eqnarray}
Here $\bar{f}_i(\x)$ denotes a degree $2\omega$ homogeneous polynomial  function  in 
$\x = (x_0,\w) \in \Real^{1+n}$ such that 
$\bar{f}_i(1,\w) = f_i(\w)$ for every $\w \in \Real^n$ $(i=0,\ldots,m)$. 

As discussed in the previous subsection, we choose an $\AC \subset \AC_{\omega}$ such that $\balpha^{\omega} \equiv (\omega,0,\ldots,0) \in \AC$ and 
the set of monomials 
$\{ \x^{\salpha+\sbeta} : \balpha \in \AC, \ \bbeta \in \AC\}$ covers all monomials involved in $\bar{f}_i(\x)$ $(i=0,\ldots,m)$, and choose 
$\Q^i \in \SymMat^{\sAC}$ $(i=0,\ldots,m)$ to satisfy
\begin{eqnarray}
\bar{f}_i(\x) = \inprod{\Q^i}{\u^{\sAC}(\x)(\u^{\sAC}(\x))^T} \  
\mbox{for every } \x \in \Real^{1+n} \ (i=0,\ldots,m). \label{eq:defQi}
\end{eqnarray}
Then, 
\begin{eqnarray}
\inprod{\Q^i}{\u^{\sAC}(\x)(\u^{\sAC}(\x))^T} = \bar{f}_i(\x) \ \mbox{for every } 
 \u^{\sAC}(\x)(\u^{\sAC}(\x))^T \in \bGamma^{\sAC} \ (i=0,\ldots,m).  
\label{eq:defQp2}
\end{eqnarray}
Define 
$ \coneJ = \left\{ \X \in \CPP^{\sAC} : \inprod{\Q^i}{\X} = 0 \ (i=1,\ldots,m) \right\}.$ 
Then, we have that 
\begin{eqnarray*}
\bGamma^{\sAC}\cap \coneJ & = & \left\{  \u^{\sAC}(\x)(\u^{\sAC}(\x))^T : \x \in \Real^{1+n}_+, \  \bar{f}_i(\x) = 0 \ (i=1,\ldots,m) \right\}. 
\end{eqnarray*}
Define 
$\H^0 \in \SymMat^{\sAC}$ such that
\begin{eqnarray*}
\H^0 & = & \mbox{the symmetric matrix in $\SymMat^{\sAC}$ whose elements are all $0$ except} \\
        &    & \mbox{the upper-leftmost element   $H^0_{\salpha^{\omega}\salpha^{\omega}}$ that is set to $1$.} 
\end{eqnarray*}
We then see that $\inprod{\H^0}{\u^{\sAC}(\x)(\u^{\sAC}(\x))^T} = x_0^{2\omega}$ for every $\x \in \Real^{1+n}$. 
It follows that 
$\X \in \bGamma^{\sAC}\cap \coneJ$ and $\inprod{\H^0}{\X} = 1$ ({\it i.e.}, $\X \in G(\bGamma^{\sAC}\cap \coneJ,1)$) 
iff there is an $\x \in \Real^{1+n}_+$ such that 
\begin{eqnarray*}
\X = \u^{\sAC}(\x)(\u^{\sAC}(\x))^T, %\ \x \in \Real^{1+n}_+, 
\ \bar{f}_i(\x)=0 \ (i=1,\ldots,m)\ \mbox{and } \ x_0^{2\omega} = 1. 
\end{eqnarray*}
Since $\x \in \Real^{1+n}_+$ implies $x_0 \geq 0$, the last equality can be replaced by $x_0 = 1$. Therefore, 
a feasible solution $\x \in \Real^{1+n}$ of POP~\eqref{eq:POP51} with the objective value $\bar{f}_0(\x)$ corresponds to 
a feasible solution $\X$ of COP($\bGamma^{\sAC}\cap \coneJ,\Q^0$) with the objective value $\inprod{\Q^0}{\X} = \bar{f}_0(\x)$ 
through the 
correspondence $\x \leftrightarrow \X = \u^{\sAC}(\x)(\u^{\sAC}(\x))^T$. 
Thus, POP~\eqref{eq:POP51} (hence POP~\eqref{eq:POP10}) is equivalent to  COP($\bGamma^{\sAC}\cap \coneJ,\Q^0)$  and $\zeta^* = \zeta(\bGamma^{\sAC}\cap \coneJ,\Q^0,1)$. 

\subsection{Reformulation of COP$(\bGamma^{\sAC}\cap\coneJ,\Q^0)$ into COP$(\coneJ,\Q^0)$}

We assume that POP~\eqref{eq:POP51} (hence \eqref{eq:POP10}) is feasible, which implies 
$G(\bGamma^{\sAC}\cap\coneJ,1) \not= \emptyset$. We also see that $\0 \not= \H^0 \in (\bGamma^{\sAC})^* \subset \coneJ^*$. 
Hence Condition I$_{\rm J}$ holds. 

Now we focus on Conditions 0$_{\rm J}$ and~II$_{\rm J}$. 
Define a sequence $\coneJ_p \subset \SymMat^{\sAC}$ $(p=0,\ldots,m)$ by~\eqref{eq:Jp40} with $\coneK=\bGamma^{\sAC}$ and 
$\coneJ_0 = {\rm co}\coneK=\CPP^{\sAC}$.
Obviously, $\coneJ = \coneJ_m$. By Lemma~\ref{lemma:face41}, we know that $\coneJ$ becomes a face of $\CPP^{\sAC}$ 
if $\Q^p$ is copositive on $\bGamma^{\sAC} \cap \coneJ_{p-1}$ $(p=1,\ldots,m)$. Thus, %we will characterize 
the copositivity of 
$\Q^p$ on $\bGamma^{\sAC} \cap \coneJ_{p-1}$  $(p=1,\ldots,m)$ can be characterized in terms of $f_i(\w)$ and 
$\bar{f}_i(\x)$ $(i=0,\ldots,m)$. 

Define 
\begin{eqnarray}
S_0 & = & \Real^{n}_+,\quad %\left\{ \w \in \Real^{n}_+ \right\}, \ 
S_p = \left\{ \w \in  S_{p-1} : f_{p}(\w) = 0 \right\} \ (p=1,\ldots,m), \label{eq:Sp} \\ 
\widetilde{S}_0 & = & \Real^{n}_+,\quad %\left\{ \w \in \Real^{n}_+ \right\}, \
\widetilde{S}_p = \left\{ \w \in  \widetilde{S}_{p-1} : \bar{f}_{p}(0,\w) = 0 \right\} \ (p=1,\ldots,m),  \label{eq:widetildeSp} \\  
\overline{S}_0 & = & \Real^{1+n}_+,\quad  %\left\{ \x \in \Real^{1+n}_+ \right\}, \
\overline{S}_p = \left\{ \x \in  \overline{S}_{p-1} : \bar{f}_{p}(\x) = 0 \right\} \ (p=1,\ldots,m). \label{eq:overlineSp}
\end{eqnarray}
By the definition of $\bGamma^{\sAC}$ and \eqref{eq:defQp2}, we observe that 
\begin{eqnarray}
\bGamma^{\sAC} \cap \coneJ_0 & = & \bGamma^{\sAC} 
= \left\{  \u^{\sAC}(\x) ( \u^{\sAC}(\x))^T: \x \in \overline{S}_0 \right\}, \nonumber \\ 
\bGamma^{\sAC} \cap \coneJ_p & = & \left\{  \u^{\sAC}(\x)( \u^{\sAC}(\x))^T  \in \coneJ_{p-1} : \bar{f}_p(\x) = 0 \right\} \nonumber  \\
& = & \left\{  \u^{\sAC}(\x) ( \u^{\sAC}(\x))^T:  \x \in \Real^{1+n}_+, \ \bar{f}_i(\x) = 0 \ (i=1,\ldots,p) \right\} \nonumber  \\
& = & \left\{  \u^{\sAC}(\x) ( \u^{\sAC}(\x))^T :  \x \in \overline{S}_p \right\} \ (p=1,\ldots,m), \label{eq:GammaJp}\\
G(\bGamma^{\sAC}\cap\coneJ_{p},\rho^{2\omega})  
& = &  \left\{  \u^{\sAC}(\x)( \u^{\sAC}(\x))^T \in \bGamma^{\sAC}\cap\coneJ_{p} :
\begin{array}{l} \x = (x_0,\w) \in \Real^{1+n}_+, \\
 x_0^{2\omega} = \rho^{2\omega}
\end{array}\right\}  \nonumber \\ 
& = &   \left\{  \u^{\sAC}(\x)( \u^{\sAC}(\x))^T :  \x = (\rho,\w)  \in \overline{S}_{p} \right\} \nonumber \\
& \ &  (\rho\geq 0, \ p=0,\ldots,m). \label{eq:Grho}
\end{eqnarray}
Now, we are ready to prove the lemma  which is used
to  establish the main theorems 
(Theorems~\ref{theorem:main50} and~\ref{theorem:main51})  with Lemma~\ref{lemma:face41} and~\ref{lemma:copositivity2}. 
\lemm \label{lemma:copositivity50} Recall that $\eta(\bGamma^{\sAC}\cap\coneJ_{p-1},\Q^p) 
= \inf \{ \langle{\Q^p},{\X}\rangle : \X\in \bGamma^{\sAC}\cap\coneJ_{p-1} \}$. We have that
\begin{description}
\item{(i) } $\eta(\bGamma^{\sAC}\cap\coneJ_{p-1},\Q^p) = \inf \left\{ \bar{f}_p(\x) : \x \in \overline{S}_{p-1} \right\}$ $ (p=1,\ldots,m)$;
\item{(ii) } $\zeta(\bGamma^{\sAC}\cap\coneJ_{p-1},\Q^p,0) = \inf \left\{ \bar{f}_p(0,\w) : \w  \in \widetilde{S}_{p-1} \right\}$ $ (p=1,\ldots,m)$;
\item{(iii) } $\zeta(\bGamma^{\sAC}\cap\coneJ_{p-1},\Q^p,1) = \inf \left\{ f_p(\w) : \w  \in S_{p-1} \right\}$ $ (p=1,\ldots,m)$;
\item{(iv) } $\zeta(\bGamma^{\sAC}\cap\coneJ_{m},\Q^0,0) = \inf \left\{ \bar{f}_0(0,\w) : \w  \in \widetilde{S}_{m} \right\}$. 
\end{description}
\elemm
\proof{
The equality in (i) follows from~\eqref{eq:defQp2}  
and~\eqref{eq:GammaJp}. 
It follows from~\eqref{eq:defQp2}, \eqref{eq:widetildeSp} 
and~\eqref{eq:Grho} with $\rho = 0$ that 
\begin{eqnarray*}
\zeta(\bGamma^{\sAC}\cap\coneJ_{p-1},\Q^p,0) & = & \inf \left\{ \inprod{\Q^p}{\X} : \X = \u^{\sAC}(\x) \u^{\sAC}(\x)^T \in G(\bGamma^{\sAC}\cap\coneJ_{p-1},0) \right\} \\ 
& = &  \inf \left\{ \bar{f}_p(\x) : \x = (0,\w)  \in \overline{S}_{p-1} \right\} \\ 
& = &  \inf \left\{ \bar{f}_p(0,\w) : \w  \in \widetilde{S}_{p-1} \right\}.  
\end{eqnarray*}
Thus we have shown (ii). For (iii), we see that 
\begin{eqnarray*}
\zeta(\bGamma^{\sAC}\cap\coneJ_{p-1},\Q^p,1) & = & \inf \left\{ \inprod{\Q^p}{\X} : \X=\u^{\sAC}(\x)(\u^{\sAC}(\x))^T \in G(\bGamma^{\sAC}\cap\coneJ_{p-1},1) \right\} \\ 
& = &  \inf \left\{ \bar{f}_p(\x) : \x = (1,\w), \w  \in \overline{S}_{p-1} \right\} \
 \mbox{(by ~\eqref{eq:defQp2}, \eqref{eq:Sp} and~\eqref{eq:Grho})} \\ 
& = &  \inf \left\{ f_p(\w) : \w  \in S_{p-1} \right\}.  
\end{eqnarray*}
(iv) follows from the same argument as the proof of (ii) with replacing $p-1$ by $m$ and $p$ by $0$. 
\qed
} 

\bigskip
We introduce the following conditions for the theorems below. 
Let $p \in \{1,\ldots,m\}$. 
\begin{eqnarray}
\inf \left\{ \bar{f}_p(\x) : \x \in \bar{S}_{p-1} \right\} & \geq & 0, \label{eq:CondPOP1-0} \\[3pt]
\inf \left\{ f_p(\w) : \w  \in S_{p-1} \right\} & \geq & 0, \label{eq:CondPOP1-1} \\ 
\inf \left\{ \bar{f}_p(0,\w) : \w  \in \widetilde{S}_{p-1} \right\} & \geq & 0, \label{eq:CondPOP1-2} \\ 
\inf \left\{ \bar{f}_0(0,\w) : \w  \in \widetilde{S}_{m} \right\} & \geq & 0. \label{eq:CondPOP2} 
\end{eqnarray}
We note that~\eqref{eq:CondPOP1-0}, \eqref{eq:CondPOP1-2} and~\eqref{eq:CondPOP2} depend on the choice of 
$\omega$, while~\eqref{eq:CondPOP1-1} is independent from the choice. 
But~\eqref{eq:CondPOP1-1} depends on  how  an optimization problem is formulated by a POP of 
the form~\eqref{eq:POP10}, as we shall  see in Section 6.2.  

\theo \label{theorem:main50} \mbox{ \ } 
Assume that $\coneJ_{p-1}$ is a face of \ $\CPP^{\sAC}$ 
for some $p\in \{1,\ldots,m\}$.
\begin{description}
\item{(i) } If~\eqref{eq:CondPOP1-0} holds, then $\coneJ_p$ is a face of  $\coneJ_{p-1}$ and a face of \ $\CPP^{\sAC}$. %$\coneJ_{p-1}$.
\item{(ii) } If~\eqref{eq:CondPOP1-1} and~\eqref{eq:CondPOP1-2} hold, then $\coneJ_p$ is a face of \ $\coneJ_{p-1}$ and a face of \ $\CPP^{\sAC}$. % $\coneJ_{p-1}$. 
\end{description}
\etheo
\proof{
(i) By (i) of Lemma~\ref {lemma:copositivity50}, we know that
$\eta(\bGamma^{\sAC}\cap\coneJ_{p-1},\Q^p) \geq 0$. Hence, the assertion follows from 
Lemma~\ref{lemma:face41}. 

(ii) Let $\coneK = \bGamma^{\sAC}$. By (ii) and (iii) of Lemma~\ref{lemma:copositivity50}, \eqref{eq:CondPOP1-1} and~\eqref{eq:CondPOP1-2}  are equivalent to 
$\zeta(\coneK\cap\coneJ_{p-1},\Q^p,1) \geq 0$ ({\it i.e}., \eqref{eq:copositive3}) and $\zeta(\coneK\cap\coneJ_{p-1},\Q^p,0) \geq 0$, 
respectively. Hence, the assertion follows from 
Lemmas~\ref{lemma:face41} and~\ref{lemma:copositivity2}.
\qed
}

\theo \label{theorem:main51} \mbox{ \ } 
Assume that POP~\eqref{eq:POP10} is feasible and that $\coneJ = \coneJ_m$ is a face of $\CPP^{\sAC}$. 
Then $\zeta^* = \zeta(\coneJ,\Q^0)$ iff 
\begin{eqnarray}
\mbox{%\eqref{eq:CondPOP2} 
$\inf \left\{ \bar{f}_0(0,\w) : \w  \in \widetilde{S}_{m} \right\} \geq 0$ {\rm ({\it i.e.}, \eqref{eq:CondPOP2})}  or $\zeta^* = -\infty$ holds. } \label{eq:iff50}
\end{eqnarray}
\etheo
\proof{
By (iv) of Lemma~\ref{lemma:copositivity50}, \eqref{eq:CondPOP2} is equivalent to $\zeta(\bGamma^{\sAC}\cap\coneJ_{m},\Q^0,0) \geq 0$, 
{\it i.e.}, 
Condition II$_{\rm J}$ with $\coneK = \bGamma^{\sAC}$. We also know that 
$\zeta^*=\zeta(\bGamma^{1+n}\cap\coneJ_m,Q^0,1)$. Hence~\eqref{eq:iff50} is equivalent to~\eqref{eq:iff31} in 
(iii) of Theorem~\ref{theorem:main}. 
\qed
}

\bigskip
Next, we make some preparations to discuss sufficient conditions for 
 \eqref{eq:CondPOP1-1}, \eqref{eq:CondPOP1-2} and~\eqref{eq:CondPOP2} to hold.
We can represent each $f_i(\w)$ as follows:
\begin{eqnarray*}
f_i(\w) & = & \hat{f}_i(\w) + \tilde{f}_i(\w), \ 
\mbox{deg$\hat{f}_i(\w) < 2\omega$ and deg$\tilde{f}_i(\w) = 2\omega$}
\end{eqnarray*}
$(i=0,1,\ldots,m)$. If deg$f_i(\w) < 2\omega$, we assume that $\tilde{f}_i(\w) \equiv 0$. 
Since $\bar{f}_i(0,\w) = \tilde{f}_i(\w)$, we know that
\begin{eqnarray*}
\bar{f}_i(0,\w) & = & \left\{\begin{array}{ll}
\mbox{a degree $2\omega$ homogeneous } \\
 \mbox{\ \hspace{8mm} polynomial function} & 
\mbox{if deg$f_i(\w)=2\omega$}, \\
0 & \mbox{otherwise, {\it i.e.}, deg$f_i(\w) < 2\omega$},
\end{array}\right.
\end{eqnarray*}
 $(i=0,\ldots,m)$. This implies that $\widetilde{S}_p$ is a cone and that $\inf \left\{ \bar{f}_0(0,\w) : \w  \in \widetilde{S}_{m} \right\} $ is either 
 $0$ or $-\infty$. Therefore, we can replace~\eqref{eq:CondPOP1-2} and~\eqref{eq:CondPOP2} by 
 \begin{eqnarray*}
 \inf \left\{ \bar{f}_p(0,\w) : \w  \in \widetilde{S}_{p-1} \right\} = 0  \ \mbox{and } 
\inf \left\{ \bar{f}_0(0,\w) : \w  \in \widetilde{S}_{m} \right\} = 0,  
 \end{eqnarray*}
respectively.

We present some sufficient conditions for~\eqref{eq:CondPOP1-1},  \eqref{eq:CondPOP1-2} and~\eqref{eq:CondPOP2} 
to hold. 
\lemm \label{lemma:barf} Let $p \in \{1,\ldots,m\}.$
\begin{description} 
\item{(i) } Assume that $f_p(\w) \geq 0$ for every $\w \in \Real^n_+$. Then~\eqref{eq:CondPOP1-1} and~\eqref{eq:CondPOP1-2} hold. 
\item{(ii) } If $\widetilde{S}_{p-1} = \{\0\}$ or deg$f_p(\w) < 2\omega$, then~\eqref{eq:CondPOP1-2} holds. 
\item{(iii) } If $\widetilde{S}_{m} = \{\0\}$ or deg$f_0(\w) < 2\omega$, then~\eqref{eq:CondPOP2} holds. 
\end{description}
\elemm
\proof{
The results in (ii) and (iii) are straightforward from the discussion above. So we only prove (i). 
Since $S_{p-1} \subset \Real^n_+$, \eqref{eq:CondPOP1-1} follows. 
To show~\eqref{eq:CondPOP1-2}, 
assume on the contrary that there is a $\tilde{\w}$ such that 
$\bar{f}_{p-1}(0,\tilde{\w}) = \tilde{f}_{p-1}(\tilde{\w}) < 0$ for some $\tilde{\w} \in \widetilde{S}_{p-1} \subset \Real^n_+$. 
Since deg$\hat{f}_{p-1}(\w) < \mbox{deg}\tilde{f}_{p-1}(\w) = 2\omega$, we have that 
\begin{eqnarray*}
\lambda\tilde{\w} \in \Real^n_+ \ \mbox{and }
f_{p-1}(\lambda\tilde{\w}) 
= \lambda^{2\omega}\left(\hat{f}_{p-1}(\lambda\tilde{\w})/\lambda^{2\omega}+ \tilde{f}_{p-1}(\tilde{\w}) \right)
< 0 
\end{eqnarray*}
for a sufficiently large $\lambda$. This contradicts the assumption. 
\qed
}

\medskip
Obviously, if $f_p(\w)$ is a sum of squares of polynomials or a polynomial with nonnegative coefficients, 
then $f_p(\w) \geq 0$ for every $\w \in \Real^n_+$. 
Otherwise,  the constraint $f_p(\w) = 0$ can be replaced by $f_p(\w)^2 = 0$ (i.e.,  the polynomial $f_p(\w)$  is  replaced by % its square
 $f_p(\w)^2$), then~\eqref{eq:CondPOP1-1} and~\eqref{eq:CondPOP1-2} are attained. 
By (ii) of Lemma~\ref{lemma:barf}, we also know that 
Condition II$_{\rm J}$ is satisfied if we take a positive integer $\omega$ such that deg$f_0(\w) < 2\omega$.  
Thus, we can easily construct an equivalent convex COP reformation, COP($\coneJ,\Q^0$) of 
POP~\eqref{eq:POP10} in theory. 

%}

%!TEX root = ./main.tex
\section{Applying Theorems~\ref{theorem:main50} and~\ref{theorem:main51} to two examples}

We illustrate how the main theorems, Theorems~\ref{theorem:main50} and~\ref{theorem:main51} 
in Section 5, can be applied to QOPs and POPs with two examples. The first one is QOP~\eqref{eq:QOP20} which has already been
reduced to COP($\bGamma^{1+n}\cap\coneJ,\Q^0$) for some cone $\coneJ \subset \CPP^{1+n} = \mbox{co}\bGamma^{1+n}$ in Section 2.2. 
The second one is a POP with some complicated %difficult 
combinatorial constraints.

%%%%%%%%%%%%%%%%%%%%%%%%%%%%%%%%%%
\subsection{QOP~\eqref{eq:QOP20} revisited}

Since  QOP~\eqref{eq:QOP20} is a special case of POP~\eqref{eq:POP10}, %we can apply 
all discussions in Section 5 can be applied to 
QOP~\eqref{eq:QOP20} if  $\u^{\sAC}(\x)$, $\bGamma^{\sAC}$ and $\CPP^{\sAC}$ are replaced
 by $(1,x_1,\ldots,x_n)$, $\bGamma^{1+n}$ and $\CPP^{1+n}$, respectively. 
%Actually
In fact, we have already constructed a cone $\coneJ \subset \CPP^{1+n} = \mbox{co}\bGamma^{1+n}$ 
and derived COP($\bGamma^{1+n}\cap\coneJ,\Q^0$), which is equivalent to QOP~\eqref{eq:QOP20}, in the same way described in 
Section 5.3. We have mentioned there that COP($\coneJ,\Q^0$) provides a CPP reformulation of QOP~\eqref{eq:QOP20} 
under conditions~\eqref{eq:binary} and~\eqref{eq:complementarity}. In this section, we prove  this fact by applying 
Theorems~\ref{theorem:main50} and~\ref{theorem:main51} . 

Recall that $\Q^p$ has been chosen to satisfy $\bar{f}_p(\x) = \inprod{\Q^p}{\x\x^T}$ for every $\x \in \Real^{1+n}$ with 
$\bar{f}_p(\x)$ given in~\eqref{eq:fbar20}  $(p=0,\ldots,m)$. 
With $\coneK = \bGamma^{1+n}$, define $\coneJ_p$, $S_p$ and $\widetilde{S}_p$ by~\eqref{eq:Jp40}, \eqref{eq:Sp} and~\eqref{eq:widetildeSp} $(p=0,\ldots,m)$, 
respectively.
Obviously $\coneJ = \coneJ_m$. We then see that 
\begin{eqnarray*}
\bar{f}_1(\x) & = & (\A\w-\b x_0)^T(\A\w-\b x_0) \geq 0 \ \mbox{for every} \ \x=(x_0,\w) \in \Real^{1+n}_+, \\
\bar{f}_2(\x) & = & \sum_{(j,k)\in I_{\rm comp}}w_jw_k \geq 0 \ \mbox{for every} \ \x=(x_0,\w) \in \Real^{1+n}_+. 
\end{eqnarray*}
By (i) of Theorem~\ref{theorem:main50}, $\coneJ_1$ and $\coneJ_2$ are faces of $\CPP^{1+n}$. Now, we show that $\coneJ_p$ is a face of $\coneJ_{p-1}$ for $p\in\{3,\ldots,m\}$. Let $p\in\{3,\ldots,m\}$ be fixed. 
It follows from~\eqref{eq:binary} % and~\eqref{eq:complementarity} 
that 
\begin{eqnarray*}
&&S_{p-1} \subset S_1 = \left\{\w \in \Real^n_+: \A\w-\b=\0 \right\} = L 
\subset \left\{\w \in \Real^n_+: w_i \leq 1 \ (i=1,\ldots,m-2)
\right\}, \\
&&\widetilde{S}_{p-1}  \subset \widetilde{S}_1 = 
\left\{\w \in \Real^n_+: \A\w=\0 \right\} = L_{\infty} \subset 
\left\{\w \in \Real^n_+:  w_i = 0 \;\; (i=1,\ldots,m-2)
\right\}. 
\end{eqnarray*}
We then see 
\begin{eqnarray*}
f_p(\w) &=& w_{p-2}(1-w_{p-2}) \geq 0 \ \mbox{for every } \w \in S_{p-1} \ \mbox{(hence \eqref{eq:CondPOP1-1} holds)}, \\
\bar{f}_p(0,\w) &=& w_{p-2}(0-w_{p-2}) = 0 \ \mbox{for every } \w \in \widetilde{S}_{p-1} \ \mbox{(hence \eqref{eq:CondPOP1-2} holds)}.
\end{eqnarray*}
By (ii) of 
Theorem~\ref{theorem:main50}, $\coneJ_p$ is a face of $\CPP^{\sAC}$. 
Thus  we have shown that $\coneJ_p$ is a face of $\CPP^{\sAC}$ for $p=3,\ldots,m$. 
Therefore, we can conclude %by (iii) of Lemma~\ref{lemma:basic} 
that $\coneJ=\coneJ_m$ is a face of $\CPP^{1+n}$. 

By Theorem~\ref{theorem:main51}, \eqref{eq:iff50} is a necessary and sufficient condition for  $\zeta^*=\zeta(\coneJ,\Q^0,1)$. 
We show that the pair of~\eqref{eq:binary} and
\eqref{eq:complementarity} is a sufficient condition for  \eqref{eq:iff50} to hold.  
We see from conditions~\eqref{eq:binary} and
~\eqref{eq:complementarity} that
\begin{eqnarray*}
\widetilde{S}_{m}  \subset \widetilde{S}_2 
&\subset& \widehat{S} \equiv \left\{\w \in \Real^n_+: 
\begin{array}{l}
\A\w = \0,\ 
w_i=0 \ (i=1,\ldots,m-2)\\
w_j=0 \ \mbox{and } w_k=0 \ ((j,k) \in I_{\rm comp})
\end{array}
\right\}.
\end{eqnarray*}
Let $\bar{\w}$ be a fasible solution of QOP~\eqref{eq:QOP20}. Suppose that
$\bar{f}_0(0,\widetilde{\w}) = \widetilde{\w}^T\C\widetilde{\w} < 0$ for some 
$\widetilde{\w} \in \widetilde{S}_{m} \red{\subset \widehat{S}}$.
%\widehat{S}$. 
Then 
$\bar{\w}+\lambda\widetilde{\w}$ is a feasible solution of 
QOP~\eqref{eq:QOP20} with the objective value $f_0(\bar{\w}+\lambda\widetilde{\w}) \rightarrow -\infty$ as $\lambda \rightarrow \infty$. Hence $\zeta^* = -\infty$. On the contrary, if there is no such a $\widetilde{\w} \in \widetilde{S}_m$, then 
$%\begin{eqnarray*}
0 %\leq \inf\{\bar{f}(0,\w) : \w \in \widehat{S} \} 
\leq \inf\{\bar{f}_0(0,\w) : \w \in \widetilde{S}_m \}; 
$ % \end{eqnarray*}
hence \eqref{eq:CondPOP2} holds. Therefore, we have shown that 
the pair of~\eqref{eq:binary}
and~\eqref{eq:complementarity} implies \eqref{eq:iff50}. 
%}

\subsection{A set of complicated combinatorial conditions from \cite{ARIMA2012}}

We consider a problem of minimizing a polynomial function in $(w_1,\ldots,w_4)$ subjct to the following combinatorial conditions. 
\begin{eqnarray}
\left. 
\begin{array}{l}
0 \leq  w_j \leq 1 \ (j=1,2,3), \ w_4 \in \{0,1\}, \\
w_1 =  1 \ \mbox{and/or } w_2  = 1,\ {\it i.e.}, \ (1-w_1)(1-w_2) = 0, \\
w_3 = 0 \ \mbox{and/or } w_1 + w_2 - w_3 = 0, \ {\it i.e.}, \ w_3(w_1 + w_2 - w_3) = 0, \\
w_4 = 0 \ \mbox{and/or } 2 - w_1 - w_2 - w_3 = 0, \ {\it i.e.}, \ w_4(2 - w_1 - w_2 - w_3) = 0.
\end{array}
\right\} \label{eq:combCondition}
\end{eqnarray}
To represent these conditions by polynomial equality constraints, we define $4$ 
polynomial functions in $\w = (w_1,\ldots,w_8) \in \Real^8$.
Choose a positive integer $\omega$ not less than the half of the degree of the objective polynomial function.
Define the $4$ polynomial functions $f_i(\w)$ in 
$\w = (w_1,\ldots,w_8) \in \Real^8$ $(i=1,\ldots,4)$ by 
\begin{eqnarray*}
& & f_1(\w) =  \sum_{k=1}^4 (w_k + w_{k+4} - 1)^{2\omega}, \
% f_2(\w) = \sum_{k=1}^4 (w_k + w_{k+4} - 1)^{2\omega} + (1-w_1)(1-w_2), \\
f_2(\w) = w_4(1-w_4) + (1-w_1)(1-w_2), \\
& & f_{3}(\w) =  w_3(w_1 + w_2 - w_3) \ \mbox{and }
f_{4}(\w) = w_4(2 - w_1 - w_2 - w_3). 
\end{eqnarray*}
Here $w_5,\ldots,w_8$ are slack variables for $w_1,\ldots,w_4$, 
respectively. %, \blue{and $\omega$ is a positive integer.}
Then the combinatorial constraints in $w_1,w_2,w_3,w_4$ above 
are satisfied iff $\w \in \Real^8_+$ and $f_i(\w) = 0$ $(i=1,2,3,4)$ for some $w_5,w_6,w_7,w_8$. 
We also set the objective polynomial function $f_0(\w)$ in 
$\w\in\Real^8$ by adding the dummy variables $w_5,w_6,w_7,w_8$ to 
the original one in $(w_1,\ldots,w_4)$. 
Thus the problem is formulated as POP~\eqref{eq:POP10} with $n=8$ and $m=4$. 
%Assuming that the degree of $f_0(\w)$ is not less than $2\omega$, 
By taking $\omega \geq \lceil \mbox{deg}f_0(\w)/2\rceil$, 
Theorems~\ref{theorem:main50} and~\ref{theorem:main51} can be applied to POP~\eqref{eq:POP10} as shown in the following.

Homogenizing the polynomial functions $f_i(\w)$ $(i=1,\ldots,4)$ 
with degree $2\omega$, we obtain that 
\begin{eqnarray*}
\bar{f}_1(x_0,\w) & =& \mbox{$\sum_{k=1}^4$}  (w_k + w_{k+4}-x_0)^{2\omega}, 
%\blue{ = \inprod{\Q^1}{\X}, }
\\[3pt]
\bar{f}_2(x_0,\w) & = & x_0^{2\omega-2}\left(w_4(x_0-w_4) + (x_0-w_1)(x_0-w_2)\right), 
% \blue{ = \inprod{\Q^2}{\X}, }
\\[3pt]
\bar{f}_{3(}x_0,\w) & = & x_0^{2\omega-2}w_3(w_1 + w_2 - w_3),  
%\blue{ = \inprod{\Q^3}{\X}, } 
\\[3pt]
 \bar{f}_{4}(x_0,\w) & = & x_0^{2\omega-2 }w_4(2x_0 - w_1 - w_2 - w_3). 
% \blue{ = \inprod{\Q^4}{\X}} 
\end{eqnarray*}
% \blue{for appropriate $\Q^i \in \SymMat^{\sAC} \;\;(i=1,\ldots,4)$ and $\X \in \bGamma^{\sAC}$.}
Then 
\begin{eqnarray*}
S_0 & = & \widetilde{S}_0 = \Real^8_+, \\ 
f_1(\w) & = & \mbox{$\sum_{k=1}^4$}  (w_k + w_{k+4}-1)^{2\omega} \geq 0 \ \mbox{for every } \w \in S_0, \\[3pt] 
S_1 & = & \left\{ \w \in \Real^8_+ : f_1(\w) = 0 \right\} \subset [0,1]^8, \\[3pt] 
f_2(\w) & = &  w_4(1-w_4) + (1-w_1)(1-w_2) \geq 0 \ \mbox{for every } \w \in [0,1]^8 \supset S_1, \\[3pt]
S_2 & = & \left\{ \w \in S_1 : f_2(\w) = 0 \right\}  = \left\{  \w \in S_1 : w_4 \in \{0,1\}, \ w_1 = 1 \ \mbox{or } w_2 = 1 \right\}, \\[3pt] 
f_{3}(\w) & = & w_3(w_1 + w_2 - w_3)  \geq 0 \ \mbox{for every } \w \in S_2, \\[3pt] 
S_3 & = & \left\{ \w \in S_2 : f_3(\w) = 0 \right\}  = \left\{  \w \in S_2 : \ w_3 = 0 \ \mbox{or } w_1 + w_2 - w_3 = 0 \right\}, \\[3pt] 
f_{4}(\w) & = & w_4(2 - w_1 - w_2 - w_3) \geq 0  \ \mbox{for every } \w \in S_3, 
\\[3pt]
\bar{f}_1(0,\w) & = & \mbox{$\sum_{k=1}^4$} (w_k + w_{k+4})^{2\omega} \geq 0 \ \mbox{for every } \w \in \widetilde{S}_0, \\[3pt] 
\widetilde{S}_1 & = & \left\{ \w \in \Real^8_+ : \bar{f}_1(0,\w) = 0 \right\} = \{\0\}, \ \widetilde{S}_2  = \widetilde{S}_3 = \widetilde{S}_4  = \{\0\}, 
\\[3pt]  
\bar{f}_p(0,\w) & \geq & 0 \ \mbox{for every } \w \in \widetilde{S}_{p-1} = \{\0\} \ (p=2,3,4), 
\\[3pt] 
\bar{f}_0(0,\w) & \geq & 0 \ \mbox{for every } \w \in \widetilde{S}_4 = \{\0\}.  
\end{eqnarray*}
Thus, we have confirmed that \eqref{eq:CondPOP1-1} and~\eqref{eq:CondPOP1-2} hold for $p=1,\ldots,4$, 
and~\eqref{eq:CondPOP2} holds with $m=4$.
By Theorems~\ref{theorem:main50} and~\ref{theorem:main51}, COP($\coneJ_4,\Q^0$)  provides a convex COP reformulation of 
POP~\eqref{eq:POP10} with $n=8$ and $m=4$. 

If additional nonnegative variables $w_9$ and $w_{10}$ are introduced and some complementarity conditions are used, then
%we  can also represent
\eqref{eq:combCondition}  can also be represented as  a single equality constraint 
\begin{eqnarray*}
\lefteqn{\sum_{k=1}^5 (w_k + w_{k+4} - 1)^{2\omega} + (w_9 -  w_1 - w_2 + w_3)^{2\omega}}  \\
& & + (w_{10} +w_1+w_2+w_3-2)^{2\omega} + 
w_4w_8 +  w_5w_6 + w_3w_9 + w_4w_{10}= 0. 
\end{eqnarray*}
In this case, we can apply the discussion at the end of Section~5.4. 

%In general,
{For given combinatorial conditions, there exist multiple ways of representing them with
binary conditions and complementarity conditions. For example,
binary conditions can be replaced by some complementarity 
conditions with slack variables. 
See  \cite{ITO2017} for more detailed discussions.

\section{Concluding remarks}

We have presented the theoretical aspects of the CPP reformulation of QOPs and its extension to POPs. To compute a lower bound
for the optimal value of POP~\eqref{eq:POP10}, numerically tractable relaxations of the problem are necessary.
Suppose that QOP~\eqref{eq:QOP20} is equivalently convexified to 
COP($\coneJ,\Q^0$) for some face $\coneJ$ of $\CPP^{1+n}$ as 
presented in Section 6.1, where 
$\coneJ$ is represented as in~\eqref{eq:coneJ10}  with $\coneK = \bGamma^{1+n}$ and co$\coneK = \CPP^{1+n}$. 
Since $\CPP^{1+n}$ is contained in the DNN cone $\SymMat^{1+n} \cap \SymN^{1+n}$, 
the CPP cone $\CPP^{1+n}$ can be relaxed to the DNN cone 
to obtain a numerically tractable DNN relaxation of QOP~\eqref{eq:QOP20} so as to 
compute a lower bound of its optimal value $\zeta_{\rm QOP}$. %of QOP~\eqref{eq:QOP20}. 
The effectiveness of this approach combined with the Lagrangian-DNN relaxation technique and 
the bisection and projection (BP) algorithm was confirmed through numerical results in 
\cite{ARIMA2017,ITO2018,KIM2013}, where 
binary QOPs, max stable set problems, multi-knapsack QOPs, quadratic assignment problems were solved. The BP algorithm was originally designed to work effectively and efficiently for 
Lagrangian-DNN relaxation problems induced from CPP reformulations of a class of QOPs with linear equality, binary and complementarity constraints in \cite{KIM2013}. 
In fact, it was shown in \cite{ARIMA2017} that Lagrangian-DNN relaxation problems induced from the CPP reformulations 
of binary QOP instances from \cite{BIQMAC} clearly provided tighter lower bounds than DNN relaxation problems  obtained from their 
standard SDP relaxations with replacing the SDP cone by the DNN cone. 

The aforementioned method using the Lagrangian-DNN relaxation technique and the BP algorithm for QOPs was 
extended to a class of sparse POPs with binary, box and 
complementarity constraints in \cite{ITO2017,KIM2016}. 
Numerical results on instances from the class showed that accurate 
lower bounds of their optimal values were efficiently obtained by the method. 
Consequently,  the theoretical study of the CPP reformulation of QOPs and its extensions to POPs  
are very important, not only for understanding of their theoretical features, but also for practical implementation. See 
\cite{ITO2017,KIM2016} for more details.

%\bibliographystyle{plain}
%\bibliography{./enhFOM}

\end{document}